\documentclass{amsart}[13pt]

\usepackage{amsrefs}
\usepackage{amsmath}
\usepackage{amssymb}
\usepackage{color}
\usepackage{hyperref}
\usepackage{stackrel}

\newtheorem{theorem}{Theorem}[section]
\newtheorem{lemma}[theorem]{Lemma}
\newtheorem{example}[theorem]{Example}
\newtheorem{definition}[theorem]{Definition}

\newtheorem{proposition}[theorem]{Proposition}
\newtheorem{remark}[theorem]{Remark}
\newtheorem{corollary}[theorem]{Corollary}

\newtheorem{nothing}[theorem]{ }

\begin{document}

\title{Induced Hausdorff metrics on quotient spaces}

\author{Ryuichi Fukuoka}
\address{Department of Mathematics, State University of Maring\'a,
87020-900, Maring\'a, PR, Brazil \\ email: rfukuoka@uem.br}
\thanks{Ryuichi Fukuoka was partially supported by the CNPq grant 305629/2012-3}

\author{Djeison Benetti}
\address{Department of Mathematics, State University of Maring\'a,
87020-900, Maring\'a, PR, Brazil \\ email: djeisonbenetti@yahoo.com.br}
\thanks{Djeison Benetti was supported by a CAPES-Funda\c c\~ao Arauc\'aria Ph.D. scholarship}

\date{September 27th, 2016}

\dedicatory{Dedicated to Professor Caio Jos\'e Colletti Negreiros on the occasion of his 60th birthday}

\begin{abstract}
Let $G$ be a group, $(M,d)$ be a metric space, $X\subset M$ be a compact subset and $\varphi:G\times M\rightarrow M$ be a left action of $G$ on $M$ by homeomorphisms. 
Denote $gp=\varphi(g,p)$. 
The isotropy subgroup of $G$ with respect to $X$ is defined by $H_X=\{g\in G; gX=X\}$. 
In this work we define the {\em induced Hausdorff metric} on $G/H_X$ by $d_X(g_1H_X,g_2H_X):=d_H(g_1X,g_2X)$, where $d_H$ is the Hausdorff distance on $M$. 
Let $\hat d_X$ be the intrinsic metric induced by $d_X$.
In this work, we study the geometry of $(G/H_X,d_X)$ and $(G/H_X,\hat d_X)$ and their relationship with $(M,d)$. 
In particular, we prove that if $G$ is a Lie group, $M$ is a differentiable manifold endowed with a metric which is locally Lipschitz equivalent to a Finsler metric, $X\subset M$ is a compact subset and $\varphi:G\times M\rightarrow M$ is a smooth left action by isometries, then $(G/H_X,\hat d_X)$ is a $C^0$-Finsler manifold. We also calculate the Finsler metric explicitly in some examples.
\end{abstract}



\maketitle

\section{Introduction}

Let $G$ be a locally compact topological group with enumerable basis and $M$ be  a Hausdorff and locally compact topological space. 
Let $\varphi: G\times M\rightarrow M$, $(g,p)\mapsto g.p:=\varphi(g,p)$, be a continuous left action of $G$ on $M$ and consider $p\in M$. 
If $\varphi$ is transitive and $H_p$ is the isotropy subgroup of $G$ with respect to $p$, then $gH_p\rightarrow gp$ is a homeomorphism from $G/H_p$ to $M$  (see \cite{Helgas}).
Moreover, if $d$ is a metric on $M$ that is compatible with its topology, then we can induce a metric $d_p$ on $G/H_p$ by $d_p(g_1H_p,g_2H_p)=d(g_1p,g_2p)$ such that $gH_P\mapsto gp$ is an isometry from $(G/H_p,d_p)$ to $(M,d)$.

In this work we study a more general situation. Consider a left action $\varphi:G\times (M,d)\rightarrow (M,d)$, where $G$ is a group, $(M,d)$ is a metric space and $\varphi$ is a left action by homeomorphisms. 
Let $X$ be a compact subset of $M$ and consider the isotropy subgroup $H_X=\{g\in G;gX=X\}$ of $G$ with respect to $X$.
Induce a metric $d_X$ on the quotient space $G/H_X$ by $d_X(g_1H_X,g_2H_X)=d_H(g_1X,g_2X)$, where $d_H$ is the Hausdorff distance on $(M,d)$ (see Proposition \ref{provainduzida}). $d_X$ is called the {\em induced Hausdorff metric} on $G/H_X$.

The intrinsic metric induced by $d_X$ (see Section \ref{preliminares}) is denoted by $\hat d_X $ and it plays an important role in this work. We are interested to study the metric spaces $(G/H_X,d_X)$ and $(G/H_X,\hat d_X)$.

Let $TM$ be the tangent bundle of a differentiable manifold $M$. 
In this work, a Finsler metric on $M$ is a continuous function $F:TM\rightarrow \mathbb R$ such that $F$ restricted to each tangent space is a norm (see \cite{Burag}, \cite{Beres1}, \cite{Beres2}, etc). 

The second important issue in this work are the metrics on differentiable manifolds that are locally Lipschitz equivalent to Finsler metrics (see Definition \ref{locallylipeqfinsl}). 
Although more general than Finsler metrics, we prove that it shares some good properties that Finsler metrics have. 
For instance all continuously differentiable curves defined on $[a,b]$ are rectifiable (Proposition \ref{locallispschitzandrectifiable}); every curve that have the same tangent vector at a point have the same speed, whenever the speed exists for one of these curves (see Definition \ref{velocidadeescalar} and Lemma \ref{mesmavelocidade}). 
We denote the space of all metrics on a differentiable manifold $M$ that are locally Lipschitz equivalent to a Finsler metric by $\mathcal L(M)$.

The main result of this work is placed in the intersection of these two subjects: 
It states that if $G$ is a Lie group, $(M,d)$ is a differentiable manifold endowed with a metric $d\in \mathcal L(M)$, $X\subset M$ is a compact subset and $\varphi:G\times (M,d)\rightarrow (M,d)$  is a smooth left action by isometries, then $\hat d_X$ is Finsler (see Theorem \ref{intrinsecaefinsler}). 
In broad sense, what happens is that the metrics $d$, $d_X$ and $\hat d_X$ are somehow related and the action is by isometries causes a regularization effect on $\hat d_X$.

In \cite{Beres1} and \cite{Beres2}, Berestovskii studies homogeneous spaces $N$ endowed with an invariant intrinsic metric $\hat d$ which induces the topology of $N$. 
One of his main results (see \cite[Theorem 3]{Beres2}) states that ``a locally compact, locally contractible homogeneous space with an intrinsic metric is isometric to a quotient space $G/H$ of some connected Lie group by a compact subgroup $H$ endowed with a Carnot-Caratheodory-Finsler metric, which is invariant with respect to the canonical action of $G$ on $G/H$.'' 
Here Carnot-Caratheodory-Finsler metric is the Finsler version of the classical Carnot-Caratheodory metric.
Moreover he proves that if $t\mapsto \exp(tv).p$ is rectifiable for every $p$ and every $v$ in the Lie algebra $\mathfrak g$ of $G$, then $\hat d$ is Finsler. 
It is worth to remark that Theorem \ref{intrinsecaefinsler} can not be settled directly from Berestovskii's Theorem because we need to prove that the topology on $G/H_X$ induced by $\hat d_X$ is the quotient topology. This result is only achieved in Theorem \ref{lipschitzlocalequivalente} at the end of this work.

This work is organized as follows.

In Section \ref{preliminares} we fix notations and give  definitions and results that are necessary for this work. The pre-requisites includes basic facts about the left action of a group on a metric space, Hausdorff distance, intrinsic metrics and Finsler metrics on differentiable manifolds. 

Let $(M,F)$ be a Finsler manifold. We define the metric $d_F: M \times M \rightarrow \mathbb R$ as 
\[
d_F(x,y)=\inf_{\gamma\in \mathcal S_{x,y}}\int F(\gamma^\prime(t))dt.
\]
where $ \mathcal S_{x,y}$ is the family of curves which connects $x$ and $y$ and are continuously differentiable by parts.
If $f:M \rightarrow M$ is a diffeomorphism, then there exist two concepts of isometry on $(M,F)$: $f:(M,d_F)\rightarrow (M,d_F)$ is an isometry of metric spaces or $df_x:T_xM \rightarrow T_{f(x)}M$ is an isometry of normed vector spaces for every $x \in M$. 
In Section \ref{isometries finsler} we prove that these two definitions coincide (Theorem \ref{nao interessa conceito de isometria}). 
In the meantime we prove that 
if $\gamma:(-\varepsilon,\varepsilon) \rightarrow (M,F)$ is path such that $\gamma^\prime(0)=v$, then the speed of $\gamma$ at $0$ is equal to $F(v)$ (See Theorem \ref{velocidade Finsler}).

In the setting of the first paragraph of this section, we have that $gH_p\mapsto gp$ is an isometry from $(G/H_X,d_X)$ to $(M,d)$. 
In more general cases, it makes sense to think that $d$, $d_X$ and $\hat d_X$ are somehow related. In Section \ref{propriedadesd} we study the influence of $\varphi$ and $d$ on $d_X$ and $\hat d_X$. Under mild conditions on $G$, $M$ and $\varphi$, we prove results such as: if $d$ and $\rho$ induces the same topology on $M$, then $d_X$ and $\rho_X$ induces the same topology on $G/H_X$ (Proposition \ref{mesmatopologiaghx}); if $d$ and $\rho$ are locally Lipschitz equivalent on $M$, then $d_X$ and $\rho_X$ (as well as $\hat d_X$ and $\hat \rho_X$) are locally Lipschitz equivalent (see Theorem \ref{lipschitzimplicalipschitz}); the quotient topology is finer than the topology induced by $d_X$ (see Proposition \ref{dxcontinua}).

In Section \ref{geometriadx} we study the geometry of $(G,d_X)$, where $G$ is a Lie group, $M$ is a differentiable manifold endowed with a metric $d$ which is compatible with the topology of $M$, $X$ is a compact subset of $M$ and $\varphi:G\times M\rightarrow M$ is a smooth left action by isometries of $G$ on $M$.
In the path to prove that $(G/H_X,\hat d_X)$ is a Finsler manifold, we need to prove that $\hat d_X$ induces the quotient topology on $G/H_X$.
In order to do that, we need to study the topology induced by $d_X$ on $G/H_X$, which can be strictly coarser than the quotient topology.
The typical example is the irrational flow on the flat torus (See Example \ref{contralocalLipschitz}).
The key result in this section is Theorem \ref{separacao bolas}, which states that there exist an $\varepsilon >0$ such that for every $gH_X \in G/H_X$, the ball $B_{d_X}(gH_X,r)=\{hH_X\in G/H_X;d_X(hH_X,gH_X)<\varepsilon\}$ is contained in a countable union of pairwise disjoint compact subsets of $G/H_X$ (compact with respect to the quotient topology). 

In Section \ref{caminhos} we study paths in $G/H_X$ endowed with the metric $d_X$ or $\hat d_X$, especially when the metric $d$ on $M$ is in $\mathcal L(M)$. 
For instance, in Proposition \ref{suaveretificavel} we prove that if $G$ is a Lie group, $(M,d)$ is a differentiable manifold endowed with a metric $d\in \mathcal L(M)$, $X\subset M$ is a compact subset and $\varphi:G\times M\rightarrow M$ is a smooth left action of $G$ on $M$, then every continuously differentiable path $\eta:[a,b]\rightarrow (G/H_X,d_X)$ is Lipschitz. 
In particular, $\eta$ is rectifiable. 
In addition, if we put the hypothesis that $\varphi$ is an action by isometries, we have that the speed of the curve $t\mapsto \exp(tv)H_X$ exists everywhere for every $v \in \mathfrak g$ (see Theorem \ref{integralkillingdiferenciavel}). 
In the way to construct the Finsler norm on $T_{H_X}G/H_X$, Lemma \ref{mesmavelocidade}, which states that curves on a differentiable manifold $(M,d\in \mathcal L(M))$ with the same derivative at some point have the same speed. Finally we prove that every path in $(G/H_X,d_X)$ is a path with respect to the quotient topology (see Theorem \ref{caminhos coincidem}). 
This result is essential to prove that the topology induced by $\hat d_X$ on $G/H_X$ is equal to the quotient topology and it uses the results that we get in Section \ref{geometriadx}.  

In Section \ref{secaofinsler} we prove Theorem \ref{intrinsecaefinsler}. 
The key result is Theorem \ref{lipschitzlocalequivalente}, that states that in the conditions of Theorem \ref{intrinsecaefinsler}, for every $gH_X\in G/H_X$, there exist a neighborhood (with respect to the quotient topology) $O\subset G/H_X$ of $gH_X$ such that $d_X\vert_{O\times O}$ and $\hat d_X\vert_{O\times O}$ are Lipschitz equivalent to a Finsler metric on $O$.
This result together with the results of Section \ref{caminhos} allow us to define the norm 
\begin{equation}
\label{velocidadeintro}
F_{H_X}(\bar v)=\lim_{\vert t\vert \rightarrow 0}\frac{d_X(c(t),H_X)}{\vert t\vert} 
\end{equation} 
on the tangent space $T_{H_X}G/H_X$, where $c:(-\varepsilon,\varepsilon)\rightarrow G/H_X$ is any curve such that $c(0)=H_X$ and $c^\prime(0)=\bar v$. 
We prove that the $G$ invariant function $F:TM\rightarrow \mathbb R$ such that $F\vert_{T_{H_X}G/H_X}=F_{H_X}$ is the Finsler metric associated to $\hat d_X$, what settles Theorem \ref{intrinsecaefinsler}.

In Section \ref{exemplos} we study some additional examples.
In particular we find an explicit expression for the Finsler metric of $(G/H_X,\hat d_X)$ when $X$ is a compact submanifold without boundary of a Riemannian manifold $M$ (see Theorem \ref{subvariedade}).

Part of this work was developed during the Ph.D. Thesis \cite{Benetti} of the second author under the supervision of the first author at State University of Maring\'a-Brazil. The authors would like to thank Anderson M. Setti and Professors Luiz A. B. San Martin, Pedro J. Catuogno, Josiney A. de Souza and Marcos R. T. Primo for their valuable suggestions. 

\section{Preliminaries}
\label{preliminares}

Classical definitions and results stated in this work can be found in \cite{Arvan},  \cite{Burag}, \cite{Carmo3}, \cite{Helgas}, \cite{Kobay1}, \cite{Kobay2}, \cite{Montg}, \cite{Munk}, \cite{Spivak1, Spivak2, Spivak3, Spivak4} and \cite{Warner}. 
{\em Every differentiable manifold will be Hausdorff with countable basis.}
  
Let $G$ be a group and consider an arbitrary non-empty set $M$. 
Denote the identity element of $G$ by $e$. 
A left action of $G$ on $M$ is a function $\varphi:G\times M \rightarrow M$, denoted by $gx:=\varphi(g,x)$, such that $ex=x$ for every $x\in M$ and $(gh)x=g(hx)$ for every $(g,h,x)\in G\times G\times M$. 
Observe that every $\varphi_g:=\varphi(g,\cdot)$ is a bijection. 
If $X\subset M$ is a subset, then the isotropy subgroup of $G$ with respect to $X$ is defined by $H_X=\{g\in G;gX=X\}$.

We can define a right action of a group on a non-empty set analogously.
{\em Unless otherwise stated, an action will stand for a left action.}

 If $(M,d)$ is a metric space with metric $d:M\times M\rightarrow \mathbb R$, then we say that $\varphi$ is an action by isometries if every $\varphi_g$ is an isometry, that is $d(x,y)=d(gx,gy)$ for every $(g,x,y)\in G\times M\times M$. Likewise we define actions by homeomorphisms and actions by diffeomorphisms.

In a metric space $(M,d)$, we denote the open ball with center $p$ and radius $r>0$ by $B_d(p,r)$, the closed ball by $B_d[p,r]$ and the sphere by $S_d(p,r)$.  
The annuli are denoted by 
\[
A_d[p,r_1,r_2]:=\{x\in M; d(p,x)\in [r_1,r_2]\}
\]
\[
A_d(p,r_1,r_2]:=\{x\in M; d(p,x)\in (r_1,r_2]\}
\]
and so on.
If $X\subset M$, then we define $B_d(X,r)=\{x\in M;d(x,X)<r\}$, $B_d[X,r]=\{x\in M;d(x,X)\leq r\}$ and $S_d(X,r)=\{x\in M;d(x,X)=r\}$, where $d(x,X)=\inf_{y\in X} d(x,y)$. 

Given a metric space $(M,d)$, the Hausdorff distance between two non-empty subsets $X_1,X_2\subset M$ is given by 
\begin{equation}
\label{defhausdorff}
d_H(X_1,X_2)=\max\{\sup_{x\in X_1}d(x,X_2),\sup_{x\in X_2}d(x,X_1)\}.
\end{equation}
The definition of $d_H$ is equivalent to
\begin{equation}
\label{hausdorffalternativo}
d_H(X_1,X_2)=\inf \{r>0;X_1\subset B_d(X_2,r) \text{ and }X_2\subset B_d(X_1,r)\}.
\end{equation}
It is straightforward to see that $d_H(X_1,X_2)=d_H(\bar X_1,\bar X_2)$, where $\bar X$ stands for the closure of $X$. 
Therefore if we want to study Hausdorff distance, we can restrict ourselves to closed subsets of $M$.

The Hausdorff distance is a metric on the set of all non-empty compact subsets of $M$ and it is an extended metric (eventually admitting infinite distance) on the set of all non-empty closed subsets of $M$.

\begin{proposition}
\label{provainduzida}
Let $G$ be a group, $(M,d)$ be a metric space, $X\subset M$ be a closed subset and $\varphi:G\times M\rightarrow M$ be an action by homeomorphisms of $G$ on $M$. Then 

\begin{enumerate}

\item $\varphi$ and $X$ induce an extended pseudometric on $G$ in the following way: 
If $g_1,g_2\in G$ and $\check d_X: G\times G\rightarrow \mathbb R$ is given by $\check d_X(g_1,g_2)=d_H(g_1X,g_2X)$, then $\check d_X$ is an extended pseudometric on $G$.
If $X$ is compact, then $\check d_X$ is a pseudometric;

\item $\check d_X$ induces an extended metric $d_X:G/H_X \times G/H_X \rightarrow \mathbb R$ given by 
\[
d_X(gH_X,hH_X):=\check d_X(g,h).
\] 
If $X$ is compact, then $d_X$ is a metric. 

\item If $\varphi$ is an action by isometries, then $d_X(agH_X,ahH_X)=d_X(gH_X,hH_X)$ for every $a\in G$ and every $gH_X,hH_X\in G/H_X$.

\end{enumerate}

\end{proposition}

{\it Proof}

\

Item 1 is straightforward from the properties of Hausdorff distance. 

Item 2 - In order to prove that $d_X$ is well defined, let $g_1H_X=g_2H_X$ and $h_1H_X=h_2H_X$. Then $g_1^{-1}g_2, h_1^{-1}h_2\in H_X$ and
\[
d_X(g_1H_X,h_1H_X)=d_H(g_1X,h_1X)
\]
\[=d_H(g_1g_1^{-1}g_2X, h_1h_1^{-1}h_2X)= d_X(g_2H_X,h_2H_X).
\]

In order to prove that $d_X$ is an extended metric, observe that the inequality $d_X\geq 0$, the symmetry of $d_X$ and the triangle inequality are straightforward. 
Moreover $d_X(gH_X,hH_X)=d_H(gX,hX)=0$ implies $gX=hX$, that is $gH_X=hH_X$. 
Therefore $d_X$ is an extended metric on $G/H_X$. Of course, if $X$ is compact, then $d_X$ is always finite and it is a metric.

Item 3 - If $\varphi$ is an action by isometries, then the equality $d_X(agH_X,ahH_X)=d_X(gH_X,hH_X)$ for every $a\in G$ and every $gH_X,hH_X\in G/H_X$ is straightforward from (\ref{defhausdorff}). $\blacksquare$

\begin{definition}
\label{definehausdorffinduzida}
Let $\varphi:G\times M\rightarrow M$ be an action by homeomorphisms of a group $G$ on a metric space $(M,d)$. 
Let $X\subset M$ be a closed subset and $H_X$ be the isotropy subgroup of $X$. 
Then the (extended) metric $d_X$ defined on $G/H_X$ in Proposition \ref{provainduzida} is called the \emph{induced Hausdorff metric} on $G/H_X$.
\end{definition}

In several situations, it is important to know when $H_X$ is a closed subset of $G$.

\begin{proposition}
\label{hxfechado}
If $\varphi:G\times M\rightarrow M$ is a continuous action of a topological group on a metric space and $X$ is closed subset of $M$, then $H_X$ is a closed subgroup of $G$. 
\end{proposition}

{\it Proof}

\

We prove that $G-H_X$ is an open subset of $G$.
The case $G=H_X$ does not need any comment. 
Consider $g\in G-H_X$. 
If there exist a $x\in X$ such that $gx\not \in X$, consider $\varphi_x:G\times \{x\}\rightarrow M$. 
Then $\varphi_x^{-1}(M-X)$ is a neighborhood of $g$ in $G-H_X$. 
Otherwise, if there is not a $x\in X$ such that $gx\not\in X$, then $gX\stackrel[\neq]{}{\subset} X$. It follows that $X\stackrel[\neq]{}{\subset} g^{-1}X$, what means that there exist a $x\in X$ such that $g^{-1}x\not\in X$. 
Using the former step there exist a neighborhood $V$ of $g^{-1}$ that does not intercept $H_X$. Thus $V^{-1}$ is a neighborhood of $g$ in $G-H_X$, what implies that $G-H_X$ is an open subset of $G/H_X$. $\blacksquare$

\

{\em Whenever we consider an (continuous or differentiable) action $\varphi:G\times (M,d)\rightarrow (M,d)$ of a Lie group (or topological group) $G$ on a differentiable manifold $M$ endowed with a metric $d$, we assume that the topology induced by $d$ will be the topology correspondent to the differentiable structure of $M$ unless otherwise stated.}

If $\varphi:G\times M\rightarrow M$ is a smooth action of a Lie group $G$ on a differentiable manifold $M$ endowed with a metric $d$ and $X\subset M$ is a closed subset, then $H_X$ is a closed subgroup of $G$ and $G/H_X$ admits an unique differentiable structure such that the natural action $\phi:G\times G/H_X\rightarrow G/H_X$ is smooth. We denote the quotient topology of $G/H_X$ by $\tau$. As a particular case, if $\varphi$ is transitive, $X=\{p\}$ and $H_p$ is the isotropy subgroup of $p$, then $gH_X\mapsto gp$ is a diffeomorphism which is also an isometry from $(G/H_X,d_X)$ to $(M,d)$. 

Let $(M,d)$ be a metric space and $\gamma:[a,b]\rightarrow M$ be a path. We denote a partition of $[a,b]$ by $\mathcal P:=\{a=t_0<t_1<\ldots<t_{n_\mathcal P}=b\}$, the norm of $\mathcal P$ by $\vert \mathcal P\vert :=\max_{i\in \{1,\ldots,n_{\mathcal P}\}}(t_i-t_{i-1})$ and
\[
\Sigma(\mathcal P):=\sum_{i=1}^{n_\mathcal P} d(\gamma(t_i),\gamma(t_{i-1})).
\]
Then length of $\gamma$ is defined by
\begin{equation}
\label{comprimento}
\ell_d(\gamma):=\sup_{\mathcal P}\Sigma(\mathcal P).
\end{equation}
If $\ell(\gamma)$ is finite, we say that $\gamma$ is rectifiable.

If $(M,d)$ is a metric space, we can define a metric $\hat d$ on $M$ by
\[
\hat d(x,y)=\inf_{\gamma \in \mathcal C_{x,y}}\ell_d (\gamma),
\]
where $\mathcal C_{x,y}$ is the family of $d$-paths connecting $x$ and $y$. 
$\hat d$ is called the intrinsic metric induced by $d$ and observe that $\hat d\geq d$ holds due to the definition of $\hat d$. 
This implies that the topology induced by $\hat d$ on $M$ is finer than the topology induced by $d$. 
In particular, every $\hat d$-path is a $d$-path.

A metric $d$ on $M$ is called {\em intrinsic} if $d=\hat d$.
We also have that $\hat d$ is always intrinsic, that is, $\hat {\hat d}=\hat d$.
If $(M,d)$ is a metric space, $\hat d$ is the intrinsic metric induced by $d$ and $\gamma:[a,b]\rightarrow M$ is a rectifiable curve in $(M,d)$, then
\begin{equation}
\label{comprimentoigual}
\ell_{\hat d}(\gamma)=\ell_d(\gamma).
\end{equation}
Moreover $\gamma$ can be reparameterized by arclength (See Proposition 2.3.12 and Proposition 2.5.9 of \cite{Burag}).

\begin{definition}
\label{velocidadeescalar}
If $(M,d)$ is a metric space and $\gamma:(a,b)\rightarrow (M,d)$ is a curve,  then the {\em speed of} $\gamma$ at $t_0\in (a,b)$ is defined by
\[
v_\gamma(t_0):=\lim_{t \rightarrow t_0}\frac{d(\gamma(t),\gamma(t_0))}{\vert t - t_0 \vert}
\]
if the limit exits.
\end{definition}

If $(M,d)$ is a metric space and $\gamma:[a,b] \rightarrow (M,d)$ is a Lipschitz curve, then the speed $v_\gamma(t)$ exists a.e. and 
\begin{equation}
\label{integral velocidade escalar}
\ell_d(\gamma)=\int_a^b v_\gamma (t) dt,
\end{equation}
where the integral above is the Lebesgue integral (See Theorem 2.7.6 of \cite{Burag}).

\begin{remark}
From now on it will be usual to endow a set with two or more different metrics and/or topologies. For instance, if a set $N$ is endowed with metrics $d_1$ and $d_2$ and a topology $\tau_N$, we use terms like $\tau_N$-open, $d_1$-bounded, $d_2$-closure, etc, in order to make clear to which topology or metric a specific term is related.
\end{remark}

Let $M$ be a differentiable manifold and denote the tangent space of $M$ at $p\in M$ by $T_pM$ and the tangent bundle of $M$ by $TM:=\{(p,v);p\in M, v\in T_pM\}$. 
Remember that a Finsler metric on $M$ is a continuous function $F:TM\rightarrow \mathbb R$ such that $F(p,\cdot):T_pM\rightarrow \mathbb R$ is a norm for every $p\in M$.
 A differentiable manifold endowed with a Finsler metric is a {\em Finsler manifold} (see \cite{Burag}, \cite{Beres1} and \cite{Beres2}). For the sake of simplicity, we write $F(v):=F(p,v)$ whenever there is no possibility of misunderstandings.

\begin{remark}
\label{outrofinsler}
There are another usual (in fact more usual) definition of Finsler manifold, where $F$ satisfies other conditions  (see, for instance, \cite{Bao} and \cite{Deng}): $F$ is smooth on $TM-TM_0$, where $TM_0=\{(p,0)\in TM;p \in M\}$, and  $F(p,\cdot)$ is a Minkowski norm on $T_pM$ for every $p\in M$. We will not use this definition in this work. Sometimes we refer to the definition we use by $C^0$-Finsler metric in order to avoid possible misunderstandings.
\end{remark}

\begin{remark}
\label{distanciafinsler}
We can define a distance function $d_F:M\times M\rightarrow \mathbb R$ from a Finsler metric by
\[
d_F(x,y)=\inf_{\gamma\in \mathcal S_{x,y}}\int F(\gamma^\prime(t))dt
\]
where $\mathcal S_{x,y}$ is the family of paths which are continuously differentiable by parts and connects $x$ and $y$.
If $(M,F)$ is a Finsler manifold and $\gamma:[a,b]\rightarrow M$ is a path which is continuously differentiable by parts, then 
\[
\ell_{d_F}(\gamma)=\int_a^b F(\gamma^\prime(t))dt,
\]
where $\ell_{d_F}(\gamma)$ is calculated according to (\ref{comprimento}) (see Subsection 2.4.2 of \cite{Burag}). Therefore 
\[
\hat d_F(x,y)
=\inf_{\gamma \in \mathcal C_{x,y}}\ell_{d_F}(\gamma)
\leq \inf_{\gamma \in \mathcal S_{x,y}}\int_a^b F(\gamma^\prime(t))dt
=d_F (x,y)
\]
and $d_F$ is intrinsic.
{\em For this reason, when a metric $d$ on a differentiable manifold is the distance function of a Finsler metric, we say by abuse of language that $d$ is Finsler.}
\end{remark}

Finally we introduce the concept of a metric which is locally Lipschitz equivalent to another metric. This concept is essential throughout this work, especially for Definition \ref{locallylipeqfinsl} below.

\begin{definition}
Let $M$ be a no-nempty set. 
Suppose that $d$ and $\rho$ are two metrics that defines the same topology on $M$.
We say that $d$ and $\rho$ are locally Lipschitz equivalent if for every $p\in M$, there exist constants $c_p,C_p>0$ and a neighborhood $V_p$ of $p\in M$ such that $c_pd(x,y)\leq \rho(x,y)\leq C_pd(x,y)$ for every $x,y\in V_p$. 
\end{definition}

\begin{definition}
\label{locallylipeqfinsl}
Let $M$ be a differentiable manifold. A metric $d:M\times M\rightarrow \mathbb R$ is {\em locally Lipschitz equivalent to a Finsler metric} if:
\begin{enumerate} 
\item the topology induced by $d$ is the topology of $M$;
\item there exist a Finsler metric $F$ on $M$ such that for every $p\in M$ there exist a neighborhood $V_p$ such that $(V_p,d_F)$ and $(V_p,d)$ are Lipschitz equivalent.
\end{enumerate}
We denote by $\mathcal L(M)$ the family of metrics on a differentiable manifold $M$ which are locally Lipschitz equivalent to a Finsler metric.
\end{definition}

The definition above does not depend on the choice of the Finsler metric because all Finsler metrics are locally Lipschitz equivalent. 

The case of a smooth action $\varphi: G \times M \rightarrow M$, where $G$ is a Lie group and $M$ is a differentiable manifold endowed with  a metric $d\in \mathcal L(M))$ is of particular importance for this work.

\section{Isometries in Finsler manifolds}

\label{isometries finsler}

When we consider a diffeomorphism $f: M\rightarrow M$ on a Finsler manifold $(M,F)$, there exist two definitions of isometry: 
In the first definition, $f: (M,d_F) \rightarrow (M,d_F)$ is an  isometry of metric spaces. The other definition states that $df_x:T_xM \rightarrow T_{f(x)}M$ is an isometry of normed vector spaces for every $x\in M$. In this section we prove that these two definitions coincide. In the meantime we prove that 
if $\gamma:(-\varepsilon,\varepsilon) \rightarrow (M,F)$ is path such that $\gamma^\prime(0)=v$, then the speed of $\gamma$ at $0$ is equal to $F(v)$ (See Theorem \ref{velocidade Finsler}).

\begin{theorem}
\label{nao interessa conceito de isometria} 
Let $f: M\rightarrow M$ be a diffeomorphism defined on a Finsler manifold $(M,F)$. 
Then $d_F(f(x),f(y))=d_F(x,y)$ for every $x,y\in M$ iff $F(f(x),df_x(v))=F(x,v)$ for every $(x,v)\in TM$.
\end{theorem}

We prove some preliminary results before the proof of Theorem \ref{nao interessa conceito de isometria}.

\

\begin{remark}
\label{considera curvas dentro}
Let $(M,d)$ be a metric space, $p\in M$ and $V$ be a $d$-neighborhood of $p$. 
Consider $r>0$ such that $B_d(p,2r)\subset V$.
Let $\gamma:[a,b]\rightarrow M$ be a rectifiable $d$-path that connects $x,y \in B_{\hat d}(p,r) \subset B_d(p,r)$ and satisfies $\gamma([a,b])\not\subset V$. Consider $\bar t\in [a,b]$ such that $\gamma(\bar t)\not\in V$. Then 
\[
\ell_{\hat d}(\gamma) =\ell_d(\gamma)\geq  \hat d(x,\gamma(\bar t))+ \hat d(y,\gamma(\bar t))\geq \hat d(p,\gamma(\bar t))- \hat d(p,x)+ \hat d(p,\gamma(\bar t)) -\hat d(p,y)
\]
\begin{equation}
\label{limitante inferior positivo}
\geq 4r - \hat d(p,x) - \hat d (p,y) \geq \hat d(x,y) + C,
\end{equation}
where $C=4r - \hat d(x,y) - \hat d(p,x) - \hat d (p,y) > 0$.
Then
\[
\hat d(x,y)
=\inf_{\gamma \in \mathcal C_{x,y}^{(M,d)}}\ell_d(\gamma)
=\inf_{\substack{\gamma \in \mathcal C_{x,y}^{(M,d)} \\ \ell_d(\gamma)
< \hat d(x,y)+C}} \ell_d(\gamma)
=\inf_{\gamma \in C_{x,y}^{(V,d)}} \ell_d(\gamma),
\]
where $\mathcal C_{x,y}^{(V,d)}$ is the family of $d$-paths that connects $x$ and $y$ and remains inside $V$. The last equality holds because we can discard the $d$-paths that does not remain in $V$ due to (\ref{limitante inferior positivo}). 
\end{remark}

\begin{remark}
\label{abertos em sequencia}
Let $\Lambda$ be a set of indices and $\{I_i\}_{i\in \Lambda}$ be a family of open intervals in $\mathbb R$ such that 
\begin{equation}
\label{nao contidos}
I_i \not\subset I_j \text{ and } I_i\not \supset I_j 
\end{equation}
or every $i\neq j$. Consider $I_i=(a_i,b_i)$ and $I_j=(a_j,b_j)$, where the cases $a_k = -\infty$ and $b_k = \infty$ are included.
Observe that if $i\neq j$ we can not have $a_i=a_j$ nor $b_i=b_j$ due to (\ref{nao contidos}) (the cases $\pm \infty$ are also included in this analysis). 

Consider, without loss of generality, that $a_i<a_j$. This forces $b_i<b_j$ due to (\ref{nao contidos}). Therefore if $i\neq j$, we can define a relationship $I_i<I_j$ if $a_i<a_j$ and $b_i<b_j$. Observe that this relationship is transitive.

\end{remark}

Consider a Finsler metric $F_0$ in $\mathbb R^n$ such that $F_0(x,v)=F_0(v)$ does not depend on $x$. In other words, $F_0$ can be identified with a norm on $\mathbb R^n$.
The proofs of Theorems \ref{nao interessa conceito de isometria} and \ref{velocidade Finsler} are based on the fact that an arbitrary Finsler metric can be locally approximated by $F_0$. 

\begin{proposition}
\label{geodesica reta}
Let $(V,F)$ be an Finsler submanifold of $(\mathbb R^n,F_0)$, that is, the inclusion map $\mathrm{inc}:V\rightarrow \mathbb R^n$ satisfies $F_0(x,d\mathrm{inc}_x(v))=F(x,v)$ for every $x\in V$ and $v\in T_xV$. Then the following statements hold:
\begin{itemize}
\item If $V$ is convex (with respect to the Euclidean metric), then 
\begin{equation}
\label{dfdf0f0}
d_F(x,y)=d_{F_0}(x,y)=F_0(x-y).
\end{equation}
In particular, straight lines in $V \subset \mathbb R^n$ are minimizers.
\item Let $(M,F)$ be a Finsler manifold which is isometric to an open subset of $(\mathbb R^n,F_0)$, that is, there exist an embedding $f:M\rightarrow \mathbb R^n$ such that $F(x,v)=F_0(f(x),df_x(v))$ for every $(x,v)\in TM$.
Then every $p \in M$ admits a neighborhood $V$ such that
\begin{equation}
\label{dfdf0f02}
d_F(x,y)=d_{F_0}(f(x),f(y))=F_0(f(x)-f(y))
\end{equation}
for every $x,y \in V$.
\end{itemize}
\end{proposition}

{\it Proof}

\

In order to prove the first item, it is enough to prove that 
\begin{equation}
\label{df0f0}
d_{F_0}(x,y)=F_0(x-y)
\end{equation}
in $(\mathbb R^n,F_0)$.
In fact, (\ref{df0f0}) implies that straight lines are minimizers in $(\mathbb R^n,F_0)$. 
Then they are minimizers also in a convex subset $V \subset (\mathbb R^n, F_0 )$ and $d_F(x,y)=F_0(x-y)$.

In order to prove the second item, it is enough to prove the first item and observe that every point of $f(M)$ admits convex neighborhoods $\tilde V$ and $\tilde W$ such that $\tilde V \subset \tilde W$ and
\[
d_F(x,y)=\inf_{\gamma \in \mathcal C_{x,y}^{\tilde W}} \ell_{d_F}(\gamma)
\]
for every $x,y \in \tilde V$ (see Remark \ref{considera curvas dentro} and notice that $d_F$ is intrinsic).

Let us prove (\ref{df0f0}).

Consider $\gamma:[0,1]\rightarrow \mathbb R^n$ given by $\gamma(t)=(1-t)x+ty$. Then $\ell_{F_0}(\gamma) = F_0(x-y)$ and $d_{F_0}(x,y)\leq F_0(x-y)$.

In order to prove that $d_{F_0}(x,y) \geq F_0(x-y)$, it is enough to prove that $\ell_{F_0}(\gamma)\geq F_0(x-y)$ for every continuously differentiable curve $\gamma:[a,b]\rightarrow \mathbb R^n$ (the continuously differentiable by parts case can be easily reduced to this case). 
Fix such a curve and consider $\varepsilon >0$. 
Then there exist a $\delta>0$ such that if $\vert t_i - t_j\vert<\delta$, then 
\[
\vert F_0(\gamma^\prime (t_i)) - F_0(\gamma^\prime (t_j)) \vert<\frac{\varepsilon}{2(b-a)}.
\] 
If we consider a partition $\mathcal P=\{a=t_0<t_1<\ldots < t_{n_\mathcal P}=b\}$ such that $\vert \mathcal P \vert<\delta$, then
\begin{equation}
\label{aproxima comprimento}
\ell_{F_0}(\gamma)=\int_a^b F_0(\gamma^\prime(t))dt=\sum_{i=1}^{n_\mathcal P}\int_{t_{i-1}}^{t_i} F_0(\gamma^\prime(t))dt\geq \sum_{i=1}^{n_\mathcal P} F_0(\gamma^\prime(\bar t_i))(t_i-t_{i-1})-\frac{\varepsilon}{2}
\end{equation}
for every $\bar t_i\in [t_{i-1},t_i]$.

For every $t\in [a,b]$, consider $\delta_t\in (0,\delta)$ such that
\begin{equation}
\label{tangentesecante}
\left\vert F_0(\gamma^\prime(t))- F_0\left( \frac{\gamma(s)-\gamma(t)}{s-t}\right)\right\vert < \frac{\varepsilon}{2(b-a)}
\end{equation}
for every $s\in [a,b]\cap (t-\delta_t,t+\delta_t)$ different from $t$. Moreover, if $t\not\in \{a,b\}$, we can choose $\delta_t$ such that $(t-\delta_t,t+\delta_t)\subset [a,b]$. 
Consider a finite subcover $\{I_j:=(\hat t_j-\delta_{\hat t_j},\hat t_j+\delta_{\hat t_j}) \}_{j=0,\ldots, k}$ of $[a,b]$, such that the intervals satisfy (\ref{nao contidos}). We can define an order on these intervals in such a way that
$I_m<I_n$ whenever $m<n$ (see Remark \ref{abertos em sequencia}). 
Observe that $\hat t_i<\hat t_j$ whenever $i<j$ because all $\hat t_k$ are in the center of the interval.
Then $\hat{\mathcal P}=\{a=\hat t_0< \ldots < \hat t_k=b \}$ is a partition of $[a,b]$ such that $I_i\cap I_{i-1}\neq \emptyset$ for every $i=1,\ldots,k$. 
Choose $s_j\in I_j\cap I_{j-1} \cap (\hat t_{j-1},\hat t_j)$. 
Then $\tilde{\mathcal P}:=\{a=\hat t_0 < s_1 < \hat t_1 < s_2 < \ldots < s_k < \hat t_k=b\}$ is a partition such that $\vert \tilde{\mathcal P} \vert<\delta$ and
\[
\ell_{F_0}(\gamma) \geq \sum_{i=1}^{n_{\mathcal P}} \left( F_0(\gamma^\prime(\hat t_{i-1}))(s_i-\hat t_{i-1}) +F_0(\gamma^\prime(\hat t_i))(\hat t_i-s_i) \right)-\frac{\varepsilon}{2}=C_1
\]
due to (\ref{aproxima comprimento}).
But (\ref{tangentesecante}) implies that
\[
F_0(\gamma^\prime(\hat t_{i-1}))( s_i - \hat t_{i-1}) >F_0\left( \gamma(s_i)-\gamma(\hat t_{i-1}) \right)-\frac{\varepsilon (s_i - \hat t_{i-1})}{2(b-a)}
\]
and
\[
F_0(\gamma^\prime(\hat t_i))( \hat t_i - s_i) >F_0\left( \gamma(\hat t_i)-\gamma(s_i) \right)-\frac{\varepsilon (\hat t_i - s_i)}{2(b-a)}.
\]
Then
\[
C_1 > \sum_{i=1}^{n_\mathcal P} \left[ F_0(\gamma(s_i)-\gamma(\hat t_{i-1})) - \frac{\varepsilon (s_i- \hat t_{i-1})}{2(b-a)} + F_0(\gamma(\hat t_i)-\gamma(s_i)) - \frac{\varepsilon (\hat t_i-s_i)}{2(b-a)} \right] -\frac{\varepsilon}{2}
\]
\[
\geq F_0(\gamma(b)-\gamma(a)) - \varepsilon=F_0(x-y) - \varepsilon.\blacksquare
\]

\begin{proposition}
\label{velocidade em fo}
Let $(M,F)$ be a Finsler manifold which is isometric (as a Finsler manifold) to an open subset of $(\mathbb R^n,F_0)$.
Consider a curve $\gamma:(-\varepsilon,\varepsilon)\rightarrow M$ such that $\gamma(0)=p$ and $\gamma^\prime(0)=v$. Then
\[
\lim_{t\rightarrow 0}\frac{d_F(\gamma(t),p)}{\vert t \vert}=F(v).
\]
\end{proposition}

{\it Proof}

\

In order to calculate
\[
\lim_{t\rightarrow 0} \frac{d_F(\gamma(t),p)}{\vert t \vert},
\]
we can restrict $\gamma$ to a neighborhood of $p$ which is isometric to a convex open subset of $(\mathbb R^n,F_0)$. 
Let $f:M\rightarrow \mathbb R^n$ be a smooth function such that $F(x,v) = F_0(f(x),df_x(v))$ for every $(x,v) \in TM$. Then

\[
\lim_{t\rightarrow 0}\frac{d_F(\gamma(t),p)}{\vert t\vert}
=\lim_{t \rightarrow 0} \frac{ F_0 (f\circ \gamma(t)-f(p))}{\vert t\vert}=C_1
\]
due to (\ref{dfdf0f02}) and
\[
C_1 
=F_0 \left( \lim_{t\rightarrow 0}\frac{f \circ \gamma(t) - f(p)}{t} \right)
=F(v)
\]
due to the continuity of $F_0$.$\blacksquare$

\begin{proposition}
\label{vizinhancas lipschitz c perto de um}
Let $F_1$ and $F_2$ be two Finsler metrics on the same underlying set $M$ and $p\in M$. Suppose that $F_1(p,\cdot)=F_2(p,\cdot)$ and fix $\varepsilon >0$. Then there exist a neighborhood $V$ of $p$ such that 
\[
(1-\varepsilon)d_{F_1}(x,y)\leq d_{F_2}(x,y)\leq (1+\varepsilon)d_{F_1}(x,y)
\]
for every $x,y\in V$.
\end{proposition}

{\it Proof}

\

If $x\in M$, then
\[
C_1(x)
=\min_{\Vert v \Vert_{F_1}=1, v\in T_xM} \frac{F_2(x,v)}{F_1(x,v)}
=\min_{v\in T_xM-\{0\}} \frac{F_2(x,v)}{F_1(x,v)}
\]
\[
C_2(x)
=\max_{\Vert v \Vert_{F_1}=1, v\in T_xM} \frac{F_2(x,v)}{F_1(x,v)}
=\max_{v\in T_xM-\{0\}} \frac{F_2(x,v)}{F_1(x,v)}
\]
are continuous functions on $M$.  
Then there exist a neighborhood $W$ of $p$ such that $C_1(x)\geq (1-\varepsilon)$ and $C_2(x)\leq (1+\varepsilon)$ for every $x\in W$, what implies that 
\begin{equation}
\label{finsler comparaveis}
(1-\varepsilon)F_1(x,v)\leq F_2(x,v) \leq (1+ \varepsilon )F_1(x,v)
\end{equation}
for every $(x,v)\in TW$.

Let $r>0$ such that $B_{d_{F_1}}(p,2r)$ and $B_{d_{F_2}}(p,2r)$ are contained in $W$. 
Remark \ref{considera curvas dentro} implies that if $x,y\in B_{d_{F_i}}(p,r)$, then
\begin{equation}
\label{basta v2}
d_{F_i}(x,y)=\inf_{\gamma \in \mathcal C_{x,y}^W} \ell_{F_i}(\gamma)
\end{equation}
for $i=1,2$, that is, in order to calculate the distance, it is enough to consider paths in $W$. But  
\begin{equation}
\label{se v2}
(1-\varepsilon)\ell_{F_2}(\gamma) \leq \ell_{F_1}(\gamma) \leq (1+\varepsilon)\ell_{F_2}(\gamma) \text{ for every }\gamma\in \mathcal C_{x,y}^{W}
\end{equation}
due to (\ref{finsler comparaveis}).
Therefore, if we set $V=B_{d_{F_1}}(p,r)\cap B_{d_{F_2}}(p,r)$, the proposition is settled.$\blacksquare$

\begin{theorem}
\label{velocidade Finsler}
Let $(M,F)$ be a Finsler manifold and $\gamma:(-a,a )\rightarrow M$ be a differentiable curve such that $\gamma(0)=p$ and $\gamma^\prime(0)=v$. Then 
\[
\lim_{t\rightarrow 0}\frac{d_F(\gamma(t),p)}{\vert t \vert}=F(v).
\]
\end{theorem}

{\it Proof}

\

Fix a neighborhood $W$ of $p$ and endow it with a metric $\tilde F$ in such a way that 
\begin{itemize}
\item $(W,\tilde F)$ is isometric to an open subset of $(\mathbb R^n,F_0)$; 
\item $\tilde F(p,\cdot)=F(p,\cdot)$.
\end{itemize}
Consider an $\varepsilon >0$. 
Proposition \ref{vizinhancas lipschitz c perto de um} states that there exist a neighborhood $V\subset W$ such that
\[
(1-\varepsilon)d_{\tilde F}(x,y)\leq d_F(x,y) \leq (1+\varepsilon) d_{\tilde F}(x,y)
\]
for every $x,y \in V$. Then
\[
 (1-\varepsilon)F(v) 
= (1-\varepsilon) \lim_{t\rightarrow 0} \frac{d_{\tilde F}(\gamma(t),p)}{\vert t \vert}
\leq \lim_{t\rightarrow 0} \frac{d_F(\gamma(t),p)}{\vert t \vert}
\leq (1+\varepsilon) \lim_{t\rightarrow 0} \frac{d_{\tilde F}(\gamma(t),p)}{\vert t \vert}
\]
\[
=(1+\varepsilon)F(v)
\]
where the first and the last equalities are due to Proposition \ref{velocidade em fo}, and the result is proved.$\blacksquare$

\

{\it Proof of Theorem \ref{nao interessa conceito de isometria}}

\

Suppose that $f:M \rightarrow M$ is a diffeomorphism such that 
\[
df_x:(T_xM,F(x,\cdot)) \rightarrow  (T_{f(x)}M,F(f(x),\cdot))
\] 
is an isometry for every $x \in M$. 
Consider $y,z\in M$. 
Then
\[
d_F(y,z)=\inf_{\gamma \in \mathcal S_{y,z}}\ell(\gamma)=\inf_{\mathcal \gamma \in S_{y,z}}\int_{a_\gamma}^{b_\gamma} F(\gamma^\prime(t)) dt
\]
\[
=\inf_{\mathcal \gamma \in S_{y,z}}\int_{a_\gamma}^{b_\gamma} F(df_{\gamma(t)}(\gamma^\prime(t))) dt
=\inf_{f\circ \gamma \in \mathcal S_{f(y),f(z)}}\ell(f\circ \gamma)=d_F(f(y),f(z)),
\]
and $f:(M,d_F) \rightarrow (M,d_F)$ is an isometry of metric spaces.

Conversely, suppose that $d_F(f(x),f(y))=d_F(x,y)$ for every $x,y\in M$. 
Consider $(x,v)\in TM$ and a differentiable curve $\gamma:(-a,a)\rightarrow M$ such that $\gamma(0)=x$ and $\gamma^\prime(0)=v$. Then
\[
F(f(x),df_x(v))
=\lim_{t \rightarrow 0}\frac{d_F(f (\gamma(t)),f(x))}{\vert t \vert}
=\lim_{t\rightarrow 0} \frac{d_F(\gamma(t),x)}{\vert t \vert}
=F(x,v),
\]
where the first and the last equality hold due to Theorem \ref{velocidade Finsler}.$\blacksquare$

\

Theorem \ref{nao interessa conceito de isometria} implies that the definition of smooth action by isometries of a Lie group on a Finsler manifold can be used without ambiguity.

\begin{corollary}
\label{nao interessa conceito de isometria 2} 
Let $\varphi:G\times M\rightarrow M$ be a smooth action of a Lie group on a Finsler manifold $(M,F)$. 
Then $d_F(gx,gy)=d_F(x,y)$ for every $(g,x,y)\in G \times M \times M$ iff $F(gx,(d\varphi_g)_x(v))=F(x,v)$ for every $(g,x,v)\in G \times M\times T_xM$.
\end{corollary}

\section{Properties of induced Hausdorff metrics}
\label{propriedadesd}

Let $G$ be a group, $(M,d)$ be a metric space, $X$ be a compact subset of $M$ and $\varphi: G \times M \rightarrow M$ be an action by homeomorphisms of $G$ on $M$. 
In this section we study the influence of $\varphi$ and $d$ on $d_X$ and $\hat d_X$. 
Under mild conditions on $G$, $M$ and $\varphi$, we prove results such as: if $d$ and $\rho$ induces the same topology on $M$, then $d_X$ and $\rho_X$ induces the same topology on $G/H_X$ (Proposition \ref{mesmatopologiaghx}); if $d$ and $\rho$ are locally Lipschitz equivalent on $M$, then $d_X$ and $\rho_X$ (as well as $\hat d_X$ and $\hat \rho_X$) are locally Lipschitz equivalent (see Theorem \ref{lipschitzimplicalipschitz}); the quotient topology is finer than the topology induced by $d_X$ on $G/H_X$ (see Proposition \ref{dxcontinua}).
These kind of results are interesting by themselves and they also help us to generalize some results that hold when $(M,F)$ is a Finsler manifold to the case $(M,d\in \mathcal L(M))$. 

\begin{proposition}
\label{dxcontinua}
Let $\varphi:G\times M\rightarrow M$ be a continuous action of a topological group $G$ on a metric space. 
Let $X\subset M$ be a compact subset. 
Then $d_X:G/H_X\times G/H_X\rightarrow \mathbb  R$ is continuous with respect to the quotient topology of $G/H_X$.
\end{proposition}

{\it Proof}

\

First of all we prove that the pseudometric $\check d_X:G\times G\rightarrow \mathbb R$ is continuous. 
Let $(g_1,g_2)\in G\times G$ and $\varepsilon>0$. We prove that there exist a neighborhood $V\times W$ of $(g_1,g_2)$ such that if $(h_1,h_2)\in V\times W$, then $\vert \check d_X(h_1,h_2) - \check d_X(g_1,g_2)\vert<\varepsilon$.

Consider the continuous function $\eta:G\times X\rightarrow \mathbb R$ given by $\eta(g,x)=d(gx,g_1x)$.
The inverse image of $(-\varepsilon/2,\varepsilon/2)$ is a neighborhood of $\{g_1\}\times X$. 
Due to the tube lemma, there exist a neighborhood $V$ of $g_1$ such that $V\times X\subset \eta^{-1}(-\varepsilon/2,\varepsilon/2)$. 
Likewise there exist a neighborhood $W$ of $g_2$ such that if $h_2\in W$, then $\vert d(h_2x,g_2x)\vert<\varepsilon /2$ for every $x\in X$.

Now observe that if $(h_1,h_2)\in V\times W$, then
\[
\check d_X(h_1,h_2)\leq \check d_X(h_1,g_1)+\check d_X(g_1,g_2) + \check d_X(g_2,h_2)
\]
\[
=d_H(h_1X,g_1X)+\check d_X(g_1,g_2)+d_H(g_2X,h_2X)
\]
\[
\leq \sup_{x\in X}d(h_1x,g_1x)+\check d_X(g_1,g_2)+\sup_{x\in X}d(h_2x,g_2x)< \check d(g_1,g_2)+\varepsilon.
\]
Analogously we have that $\check d_X(g_1,g_2)< \check d_X(h_1,h_2)+\varepsilon$ and we have proved that for every $(h_1,h_2)\in V\times W$ we have that $\vert \check d_X(h_1,h_2)-\check d_X(g_1,g_2)\vert<\varepsilon$. 
Thus $\check d_X$ is continuous.

Finally let $(g_1H_X, g_2H_X)\in G/H_X\times G/H_X$ and fix $\varepsilon >0$. 
Consider a neighborhood $V\times W$ of $(g_1,g_2)\in G\times G$. 
Observe that $\pi(V)\times \pi(W)$ is a neighborhood of $(g_1H_X, g_2H_X)\in G/H_X \times G/H_X$ such that if $(h_1H_X,h_2H_X) \in \pi(V)\times \pi(W)$, then $\vert d_X (h_1H_X,h_2H_X)-d_X(g_1H_X,g_2H_X)\vert = \vert \check d_X(h_1,h_2)- \check d_X(g_1,g_2)\vert <\varepsilon$.
Therefore $d_X$ is continuous with respect to the quotient topology on $G/H_X$. $\blacksquare$

\

\begin{remark}
\label{curvanoquocientecontinua}
If $(M,\tau_M)$ is a topological space and $d$ is a metric which is not necessarily compatible with $\tau_M$, then 
$d:M\times M\rightarrow \mathbb R$ is $\tau_M$-continuous if and only if $\tau_M$ is finer than the topology induced by $d$ (see \cite{Munk}). 
In Proposition \ref{dxcontinua}, $d_X$ is continuous with respect to the quotient topology $\tau$ of $G/H_X$.
Therefore the quotient topology is finer than the topology induced by $d_X$. In particular, if $\alpha$ is a path in $(G/H_X,\tau)$, then $\alpha$ will be also a path in the metric space $(G/H_X,d_X)$. 
\end{remark}

\begin{example}
\label{intrinsecanaocontinua}
In the conditions of Proposition \ref{dxcontinua}, $\hat d_X$ is not necessarily continuous with respect to the quotient metric on $G/H_X$. 
As an example, consider $M=\mathbb R$ with the metric $d(x,y)=\sqrt{\vert x-y\vert}$.
It is not difficult to show that $\hat d(x,y)=\infty$ if $x\neq y$.
If we consider the natural action of the group $G=(\mathbb R,+)$ on $M$ by addition and we take $X=\{ 0 \}$, then $H_X=\{0\}$ and $(G/H_X,d_X)$ is isometric to the metric space $M$. 
Then $(G/H_X,\hat d_X)$ is isometric to $(M,\hat d)$ and $\hat d_X$ is not continuous with respect to the quotient topology on $G/H_X$.
\end{example}

The next step is to prove Proposition \ref{mesmatopologiaghx} that states that under mild conditions, if $d$ and $\rho$ are metrics that induces the same topology on $M$, then $d_X$ and $\rho_X$ induces the same topology on $G/H_X$. We will prove some lemmas before.

\begin{lemma}
\label{brcompacto}
Let $(M,d)$ be a locally compact metric space and $Y$ be a compact subset of $M$. Then there exist a $r>0$ such that $\bar B_d(Y,r)$ is compact.
\end{lemma}

{\it Proof}

\

Let $\{U_i\}_{i=1,\ldots k}$ be a finite open covering of $Y$ such that $\bar U_i$ is compact for every $i=1,\ldots,k$. Then $A:=\cup_{i=1}^k U_i$ is an open set containing $Y$ such that $\bar A=\cup_{i=1}^k \bar U_i$ is compact. If $A=M=\bar A$, then there is nothing to prove. 
If $A\neq M$, then there exist a $r>0$ such that $B_d(Y,r)\subset A$ because $Y$ is compact and $d(\cdot,M-A):Y\rightarrow \mathbb R$ is continuous. Therefore $\bar B_d(Y,r)\subset \bar A$ is compact.$\blacksquare$

\

\begin{lemma}
\label{compactoepsilondeltauniforme} Let $d$ and $\rho$ be two metrics that induces the same topology on $M$. Let $Y$ be a compact subset of $M$. Then for every $\varepsilon>0$, there exist a $\delta>0$ such that $B_\rho(x,\delta)\subset B_d(x,\varepsilon)$ for every $x\in Y$. In particular, $B_\rho(Z,\delta)\subset B_d(Z,\varepsilon)$ for every $Z\subset Y$. 
\end{lemma}

{\it Proof}

\

Fix $\varepsilon > 0$. For every $x\in Y$, consider $\delta_x>0$ such that $B_\rho(x,\delta_x)\subset B_d(x,\varepsilon/2)$. 
Let $2\delta$ be the Lebesgue number of the open cover $\{B_\rho(x,\delta_x)\}_{x\in Y}$ of $Y$. 
Then $B_\rho(x,\delta)\subset B_\rho(z,\delta_z)$ for some $z\in Y$ and $B_\rho(z,\delta_z)\subset B_d(z,\varepsilon/2)\subset B_d(x,\varepsilon)$. The last statement is straightforward. $\blacksquare$

\

Of course we can change the roles of $\rho$ and $d$ in Lemma \ref{compactoepsilondeltauniforme}.

\begin{proposition}
\label{mesmatopologiaghx}
Let $d$ and $\rho$ be metrics on $M$ such that they induces the same topology $\tau_M$ on $M$ and suppose that $(M,\tau_M)$ is locally compact.
Consider an action $\varphi:G\times M\rightarrow M$ by homeomorphisms of a topological group $G$ on $M$. 
Let $X$ be a compact subset of $M$. Then $d_X$ and $\rho_X$ induces the same topology on $G/H_X$.  
\end{proposition}

{\it Proof}

\

Fix $gH_X\in G/H_X$ and consider $r >0$ such that $\bar B_d(gX,r)$ is compact (See Lemma \ref{brcompacto}). 
Let $\varepsilon>0$ and consider $\delta>0$ such that $B_\rho(x,\delta)\subset B_d(x,\varepsilon/2)$ for every $x\in \bar B_d(gX,r)$ (see Lemma \ref{compactoepsilondeltauniforme}). 
Without loss of generality, we can consider $\varepsilon<r$.
If $hH_X \in B_{\rho_X}(gH_X,\delta)$, then $\rho_H(gX,hX)<\delta$ what implies that $gX\subset B_\rho(hX,\delta)$ and $hX\subset B_\rho(gX,\delta)$ (see (\ref{hausdorffalternativo})). 
Consequently we have that $gX\subset B_d(hX,\varepsilon/2)$ and $hX\subset B_d(gX,\varepsilon/2)$ due to the definition of  $\delta$ and we have that $hH_X\in B_{d_X}(gH_X,\varepsilon)$ due to (\ref{hausdorffalternativo}). 
Therefore given $gH_X\in G/H_X$ and $\varepsilon>0$, there exist $\delta>0$ such that $B_{\rho_X}(gH_X,\delta)\subset B_{d_X}(gH_X,\varepsilon)$. 
Thus the topology induced by $\rho_X$ is finer than the topology induced by $d_X$. 
If we change the roles of $\rho_X$ and $d_X$, we get that the topology induced by $d_X$ is finer than the topology induced by $\rho_X$.$\blacksquare$

\

The next step is to prove Theorem \ref{lipschitzimplicalipschitz}. We prove some preliminary results first.

\begin{proposition}
\label{desigualdadeinduzigualchapeu}
Let $d$ and $\rho$ be two metrics defined on $M$ which are locally Lipschitz equivalent. 
Then $\hat d$ and $\hat \rho$ are locally Lipschitz equivalent. 
\end{proposition}

{\it Proof}

\

We denote the topology induced by $d$ (or $\rho$) on $M$ by $\tau_M$. 
Let $p\in M$, $V$ be a $\tau_M$-neighborhood of $p$ such that $d\vert_V$ and $\rho\vert_V$ are Lipschitz equivalent. 
Let $c,C>0$ such that 
\[
cd(x,y) \leq \rho(x,y) \leq Cd(x,y)
\]
for every $x,y\in V$. 
Notice that
\begin{equation}
\label{compara comprimento}
c\ell_d(\gamma) \leq \ell_\rho (\gamma) \leq C \ell_d (\gamma)
\end{equation}
for every $\tau_M$-path $\gamma$ in $V$.

Consider $r>0$ such that $B_d(p, 2 r)\subset V$ and $B_\rho(p,2r)\subset V$. 
Remember that $B_{\hat d}(p,r)\subset B_d(p, r)$ and $B_{\hat \rho}(p,r)\subset B_\rho(p,r)$.

Due to Remark \ref{considera curvas dentro}, if $x,y\in B_{\hat d}(p,r)$, then
\begin{equation}
\label{d chapeu}
\hat d(x,y)
=\inf_{\gamma \in \mathcal C^{(M,\tau_M)}_{x,y}} \ell_d(\gamma)
=\inf_{\gamma \in \mathcal C^{(V,\tau_M)}_{x,y}} \ell_d(\gamma),
\end{equation}
where $\mathcal C^{(V,\tau_M)}_{x,y}$ is the family of $\tau_M$-paths in $V$ that connects $x$ and $y$. 
Analogously we have that
\begin{equation}
\label{rho chapeu}
\hat \rho(x,y)
=\inf_{\gamma \in \mathcal C^{(M, \tau_M)}_{x,y}} \ell_\rho(\gamma)
=\inf_{\gamma \in \mathcal C^{(V,\tau_M)}_{x,y}} \ell_\rho(\gamma)
\end{equation}
for every $x,y\in B_{\hat \rho}(p,r)$.

Now we will prove that for every $\varepsilon$, there exist a $\delta$ such that $B_{\hat \rho}(p,\delta) \subset B_{\hat d}(p,\varepsilon)$. Without loss of generality we can consider $\varepsilon\leq r$. If $\delta=c \varepsilon$ and $x\in B_{\hat \rho}(p,\delta)$, then
\[
\hat d(p,x)
= \inf_{\gamma \in \mathcal C^{(V,\tau_M)}_{p,x}} \ell_d(\gamma)
\leq \frac{\inf_{\gamma \in \mathcal C^{(V,\tau_M)}_{p,x}} \ell_\rho(\gamma) }{c}
=  \frac{\hat \rho (p,x)}{c} 
< \varepsilon
\]
due to (\ref{compara comprimento}), and $x\in B_{\hat d}(p,\varepsilon)$. 
Therefore $B_{\hat \rho}(p,\delta) \subset B_{\hat d}(p,\varepsilon)$ and $\hat \rho$ is finer than $\hat d$.
The proof that $\hat d$ is finer than $\hat \rho$ is analogous. Then the topologies induced by $\hat d$ and $\hat \rho$ coincide.

Finally in order to see that $\hat d$ and $\hat \rho$ are Lipschitz equivalent on $W=B_{\hat d}(p, r)\cap B_{\hat \rho}(p,r)$, it is enough to observe that the relationships
\[
c\hat d(x,y) \leq \hat \rho(x,y) \leq C \hat d(x,y)
\]
hold for every $x,y\in W$ due to (\ref{compara comprimento}), (\ref{d chapeu}) and (\ref{rho chapeu}).$\blacksquare$

\

\begin{proposition}
\label{compactolipschitzequivalente}
Let $d$ and $\rho$ be two locally Lipschitz equivalent metrics on $M$ and consider a compact subset $Y\subset M$. Then $d\vert_{Y\times Y}$ and $\rho\vert_{Y\times Y}$ are Lipschitz equivalent.
\end{proposition}

{\it Proof}

\

Consider an open cover $\{V_i\}_{i=1,\ldots k}$ of $Y$ and positive constants $\bar c,\bar C>0$ such that 
\begin{equation}
\label{lipschitzlocal}
\bar c.d(x,y)\leq \rho(x,y)\leq \bar C.d(x,y) \text{ for every } x,y\in V_i,
\end{equation}
where (\ref{lipschitzlocal}) holds for every $i=1,\ldots, k$.

$A:=Y\times Y-\cup_{i=1}^k V_i \times V_i$ is a compact subset of $M\times M$ and $d\vert_A$ as well as $\rho\vert_A$ admit strictly positive lower bounds. 
Then the quotients $(d/\rho)\vert_A$ and $(\rho/ d)\vert_A$ assume their supremum and infimum and we have that 
\begin{equation}
\label{complemento}
\tilde c d(x,y)\leq \rho(x,y)\leq \tilde C d(x,y) \text{ for every } (x,y)\in A
\end{equation}
for some positive constants $\tilde c$ and $\tilde C$. From (\ref{lipschitzlocal}) and (\ref{complemento}) we have that there exist constants $c,C>0$ such that
\begin{equation}
\label{lipschitzcompacto}
c d(x,y)\leq \rho(x,y)\leq C d(x,y) \text{ for every } x,y\in Y.\blacksquare
\end{equation}

\

\

\

\begin{proposition}
\label{hausdorfflipschitz}

\

\begin{enumerate}
\item Let $d$ and $\rho$ be two metrics on $M$ such that $d(x,y)\leq C\rho(x,y)$ for some $C>0$ and every $x,y\in M$. 
If $X,Y\subset M$, then $d_H(X,Y)\leq C\rho_H(X,Y)$.
\item Let $G$ be a topological group, $d$ and $\rho$ be two Lipschitz equivalent metrics on $M$, $X\subset M$ be a compact subset and $\varphi:G\times M\rightarrow M$ be a continuous action of $G$ on $M$.
Then $d_X$ and $\rho_X$ are also Lipschitz equivalent.
\end{enumerate}
\end{proposition}

{\it Proof}

\

Item 2 is a direct consequence of Item 1. 

We prove that if $d(x,y)\leq C\rho(x,y)$ for every $x,y\in M$ and $X,Y$ are non-empty subsets of $M$, then

\[
d_H(X,Y)\leq C\rho_H(X,Y).
\]
Remember that
\[
d_H(X,Y)=\max\{\sup_{x\in X}\inf_{y\in Y}d(x,y),\sup_{y\in Y}\inf_{x\in X}d(x,y)\}.
\]
Then
\[
d(x,y)\leq C\rho(x,y)\text{ for every }x\in X,y\in Y;
\]
\[
\inf_{\bar y\in Y}d(x,\bar y)\leq C \rho(x,y)\text{ for every }x\in X,y\in Y;
\]
\[
\inf_{y\in Y}d(x,y)\leq C\inf_{y\in Y}\rho(x,y)\text{ for every }x\in X;
\]
\[
\inf_{y\in Y}d(x,y)\leq C\sup_{\bar x\in X}\inf_{y\in Y}\rho(\bar x,y)\text{ for every }x\in X;
\]
\[
A_d:=\sup_{x\in X}\inf_{y\in Y}d(x,y)\leq C\sup_{x\in X}\inf_{y\in Y}\rho(x,y):=CA_\rho.
\]
Analogously we have that
\[
B_d:=\sup_{y\in Y}\inf_{x\in X}d(x,y)\leq C\sup_{y\in Y}\inf_{x\in X}\rho(x,y):=CB_\rho.
\]
Finally
\[
A_d,B_d\leq \max\{CA_\rho,CB_\rho\}
\]
and
\[
\max\{A_d,B_d\}\leq C\max\{A_\rho,B_\rho\}.\blacksquare
\]

\

\begin{theorem}
\label{lipschitzimplicalipschitz}
Let $d$ and $\rho$ be two locally Lipschitz equivalent metrics on $M$ and suppose that $M$ is locally compact (with respect to $d$ or $\rho$). 
Let $\varphi:G\times M\rightarrow M$ be a continuous action of a topological group $G$ on $M$ and consider a compact subset $X\subset M$. 
Then
\begin{enumerate} 
\item $d_X$ and $\rho_X$ are locally Lipschitz equivalent on $G/H_X$;
\item $\hat d_X$ and $\hat \rho_X$ are locally Lipschitz equivalent on $G/H_X$.
\end{enumerate}
\end{theorem}

{\it Proof}

\

The second item follows from the first item and Proposition \ref{desigualdadeinduzigualchapeu}. 
Then it is enough to prove the first item.

First of all $d_X$ and $\rho_X$ induces the same topology on $G/H_X$ due to Proposition \ref{mesmatopologiaghx}.

For $gH_X\in G/H_X$, consider $r>0$ such that $\bar B_d(gX,r)$ is compact (see Lemma \ref{brcompacto}). We claim that $d_X$ and $\rho_X$ are Lipschitz equivalent on $B_{d_X}(gH_X,r)$. 
In fact, notice that $d$ and $\rho$ are Lipschitz equivalent on $A=B_d(gX,r)$ (Proposition \ref{compactolipschitzequivalente}). Then
\begin{equation}
\label{ctelipschitz}
cd(x,y)\leq \rho(x,y)\leq Cd(x,y)
\end{equation} 
for every $x,y\in A$ and some positive constants $c,C>0$. 
If we consider $h_1H_X,h_2H_X\in B_{d_X}(gH_X,r)$, then $h_1X,h_2X\subset B_d(gX,r)$ (see (\ref{hausdorffalternativo})) and 
\[
cd_X(h_1H_X,h_2H_X) =cd_H(h_1X,h_2X)\leq \rho_H(h_1X,h_2X)= \rho_X (h_1H_X,h_2H_X)
\]
\[
\leq Cd_H(h_1X,h_2X)=Cd_X(h_1H_X,h_2H_X)
\]
due to Proposition \ref{hausdorfflipschitz}.$\blacksquare$

\section{Geometry and topology of $(G/H_X,d_X)$}
\label{geometriadx}

Let $G$ be a Lie group, $M$ be a differentiable manifold endowed with a metric $d$ and $X$ be a compact subset of $M$. 
Consider an action $\varphi:G\times M\rightarrow M$ by isometries.
In this section, we study the geometry and topology of $(G/H_X,d_X)$. We know that $d_X$-open subsets of $G/H_X$ are open in the quotient topology $\tau$ (Proposition \ref{dxcontinua}). 
In this section we prove that there exist an $\varepsilon >0$ such that $B_{d_X}(gH_X,\varepsilon)$ is contained in a countable family of pairwise disjoint $\tau$-compact subsets of $G/H_X$ (Theorem \ref{separacao bolas}).
This result will be important in order to prove that paths in $(G/H_X,d_X)$ are paths in $(G/H_X,\tau)$ (Theorem \ref{caminhos coincidem}), which is used in order to prove the second item of Theorem \ref{lipschitzlocalequivalente}. More precisely, if $\tau_X$ is the topology of $(G/H_X , d_X)$ and $\hat\tau_X$ is the topology of  $(G/H_X, \hat d_X)$, then $\tau_X\subset \hat \tau_X =\tau$.

The following example illustrates Theorem \ref{separacao bolas}.

\begin{example}
\label{contralocalLipschitz}
Consider the Lie group $G=(\mathbb R,+)$ with the canonical differentiable structure, and let $M=(\mathbb R\times \mathbb R)/(\mathbb Z\times \mathbb Z)$ be the flat torus. 
Denote its metric by $d$.
We represent a point in $M$ by $(\bar x,\bar y)$, where $x,y\in \mathbb R$ and $\bar x$ is the equivalence class of $x\in \mathbb R$ in $\mathbb R/\mathbb Z$.
Consider the action $\varphi:G\times M\rightarrow M$ given by $t(\bar x,\bar y)=(\overline{t+x},\overline{t\sqrt 2+y})$, which is a irrational flow on the flat torus.
Consider $X=\{(\bar 0,\bar 0)\}$. Then $H_X=\{0\}$ and $G/H_X\cong G$. 
Given an $\varepsilon>0$, the open ball $B_{d_X}(0,\varepsilon)$ is not bounded with respect to the Euclidean metric on $G$. 
In fact, if we take an arbitrarily big $N>0$, there exist a $t>N$ such that $t(\bar 0,\bar 0)\in B_d((\bar 0,\bar 0),\varepsilon)$. Then $t\in B_{d_X}(0,\varepsilon)$ and $B_{d_X}(0,\varepsilon)\subset \mathbb R$ is unbounded with respect to the canonical metric for every $\varepsilon>0$. 
Therefore $\tau_X\underset{\neq}{\subset}\tau$.
It is also easy to notice that for a sufficiently small $\varepsilon$, $B_{d_X}(0,\varepsilon)$ has infinite countable path-connected components.
\end{example}

Before the proof of Theorem \ref{separacao bolas}, we present a general setting in order to be referenced afterwards because it will be used frequently.

\begin{nothing}
\label{mhx}
{\rm\bf General setting}
\end{nothing}

Let $G$ be a Lie group, $(M,d)$ be a differentiable manifold endowed with a metric $d$, $X$ be a compact subset of $M$ and $\varphi:G\times M\rightarrow M$ be a smooth action. 
Consider a decomposition $\mathfrak g=\mathfrak h_X \oplus \mathfrak m$ of the Lie algebra $\mathfrak g$ of $G$, where $\mathfrak h_X$ is the Lie algebra of the isotropy subgroup $H_X$ and $\mathfrak m$ is a subspace of $\mathfrak g$.

Fix an Euclidean metric $d_\mathfrak m$ on $\mathfrak m$.
Let $r>0$ such that $\pi \circ \exp$ restricted to $B_{d_\mathfrak m}(0,r)\subset \mathfrak m$ is a diffeomorphism over its image.

The map $\pi\circ \exp$ induce a metric on $B_{d_\mathfrak m}(H_X,r):=\pi\circ \exp (B_{d_\mathfrak m}(0,r))$ which for the sake of simplicity we call $d_\mathfrak m$. 
We will consider three different metrics on $B_{d_\mathfrak m}(H_X,r)$: $d_\mathfrak m$, $d_X$ and $\hat d_X$. $\blacksquare$

\begin{theorem}
\label{separacao bolas}
Let $G$ be a Lie group, $(M,d)$ be a differentiable manifold endowed with a metric $d$, $X$ be a compact subset of $M$ and $\varphi:G\times M\rightarrow M$ be a smooth action by isometries. 
Then there exist an $\varepsilon >0$ such that every ball $B_{d_X}(gH_X,\varepsilon)$ is contained in a countable (finite or infinite) union $\cup_{i\in \Lambda}K_i$ of pairwise disjoint $\tau$-compact subsets.
\end{theorem}

{\it Proof}

\

Due to the $G$-invariance of $d_X$ on $G/H_X$, it is enough to prove that there exist an $\varepsilon >0$ such that $B_{d_X}(H_X,\varepsilon)$ is contained in a countable union $\cup_{l\in \Lambda}K_l$ of disjoint $\tau$-compact subsets.

Suppose that we are in the settings of $\ref{mhx}$.
For every $g\in G$, define $B^g_{d_\mathfrak m}[gH_X,s_1]:=gB_{d_{\mathfrak m}}[H_X,s_1]$ and $B^g_{d_\mathfrak m}(gH_X,s_1):=g B_{d_{\mathfrak m}}(H_X,s_1)$, where $s_1\in (0,r)$. Analogously we define
\[
A_{d_\mathfrak m}^g[gH_X,s_1, s_2]=gA_{d_\mathfrak m}[H_X,s_1,s_2],
\]
\[
A^g_{d_\mathfrak m}[gH_X,s_1, s_2)=gA_{d_\mathfrak m}[H_X,s_1,s_2),
\]
and so on.  

Fix $s\in (0,r)$ and set
\[
A_1^g=A_{d_\mathfrak m}^g[gH_X,s/2,s],
\]
\[
R_1^g=\min_{a\in A_1^g}d_X(a,gH_X).
\]
Consider $\tilde s\in (0,s)$ such that
\[
R_2^g=\max_{a\in B_{d_\mathfrak m}^g[gH_X,\tilde s]}d_X(a,gH_X)\in \left( 0, \frac{R_1^g}{2}\right).
\]
For the sake of simplicity, denote $B_2^g:=B_{d_\mathfrak m}^g[gH_X,\tilde s]$.
Finally set
\[
A_3^g=A_{d_\mathfrak m}^g[gH_X,\tilde s/2, s]
\]
and
\[
R_3^g=\min_{a\in A_3^g}d_X(a,gH_X).
\]

If $g\in G$, then
\[
R_1^g:=\min_{\hat a\in A_1^g}d_X(\hat a,gH_X)
=\min_{a\in A_1^e}d_X(ga,gH_X)
\]
\[
=\min_{a\in A_1^e}d_X(a,H_X)=R_1^e
\]
and analogously $R_2^g=R_2^e$ and $R_3^g=R_3^e$ for every $g\in G$. Denote these constants by $R_1$, $R_2$ and $R_3$ respectively.

All these subsets are defined in order to find an $\varepsilon >0$ such that $B_{d_X}(H_X,\varepsilon)$ is contained in a pairwise disjoint union $\cup_{i\in \Lambda} K_i$ of $\tau$-compact subsets. 
More precisely, we prove that everything works if we choose $\varepsilon=\min\{R_2/2,R_3/2\}$ and $K_i=B_2^{g_i}$, where $g_i\in B_{d_X}(H_X, \varepsilon)$ are chosen strategically.

Observe that $G/H_X-A_1^g$ can be written as the disjoint union $B_{d_\mathfrak m}^g(gH_X,s/2)\cup (G/H_X-B_{d_\mathfrak m}^g[gH_X,s])$, that is, one inside part and one outside part. 
The same type of decomposition holds for 
\[
G/H_X-A_3^g=B_{d_\mathfrak m}^g(gH_X,\tilde s/2)\cup (G/H_X-B_{d_\mathfrak m}^g[gH_X,s]).
\] 
Moreover $B_{d_X}(gH_X,R_3)\cap A_3^g=\emptyset$ for every $g\in G$ due to the definition of $R_3$, what implies that if $O$ is a $\tau$-connected component of $B_{d_X}(gH_X,R_3)$, then either $O\subset B_{d_\mathfrak m}^g(gH_X,\tilde s/2)$ (inside part) or else $O\subset G/H_X-B_{d_\mathfrak m}^g[gH_X,s]$ (outside part). 

It is clear that $A_1^g \cap B_2^g = \emptyset$ due to the definition of $B_2^g$. We claim that 
\begin{equation}
\label{intersecao vazia a b}
A_1^{g_1}\cap B_2^{g_2}=\emptyset
\end{equation}
whenever $g_1H_X$ and $g_2H_X$ are close enough. More precisely, (\ref{intersecao vazia a b}) holds if  $g_1H_X,$ $g_2H_X\in B_{d_X}(H_X, R_2/2)$. 
In fact, observe that if $\hat a_1\in A_1^{g_1}$ and $\hat a_2\in B_2^{g_2}$, then 
\[
d_X(\hat a_2,\hat a_1 )=d_X(g_2a_2,g_1a_1)=C_1
\]
for some $a_1\in A_1^e$ and $a_2\in B_2^e$. Consequently
\[
C_1\geq d_X(g_1H_X,g_1a_1)-d_X(g_1H_X,g_2H_X)-d_X(g_2H_X,g_2a_2)
\]
\[
= d_X(H_X,a_1)-d_X(g_1H_X,g_2H_X)-d_X(H_X,a_2) \geq R_1-2R_2>0
\]
due to the $G$-invariance of $d_X$, what proves the claim. 
Therefore the (connected) subset $B_2^{g_2}$ is contained either in $B_{d_\mathfrak m}(g_1H_X,s/2)$ or else in $G/H_X-B_{d_\mathfrak m}[g_1H_X,s]$. 
In particular, for every $g_1H_X,g_2H_X \in B_{d_X}(H_X,R_2/2)$ we have that
\begin{equation}
\label{bi opcoes} 
B_2^{g_2}\subset B_{d_\mathfrak m}^{g_1}[g_1H_X,s] \text{ or else } B_2^{g_2}\cap B_{d_\mathfrak m}^{g_1}[g_1H_X,s]=\emptyset,
\end{equation}
what implies that the $\tau$-compact subsets $B_2^{g_2}$ that are contained in $B^{g_1}_{d_\mathfrak m}[g_1H_X,s]$ does not intercept those $B_2^{g_2^\prime}$ that are not contained in $B^{g_1}_{d_\mathfrak m}[g_1H_X,s]$.

Set $\varepsilon=\min\{R_3/2,R_2/2\}$. Notice that for every $gH_X\in B_{d_X}(H_X,\varepsilon)$, we have that 
\[
A_3^g\cap B_{d_X}(H_X,\varepsilon)=\emptyset
\] 
because $B_{d_X}(H_X,\varepsilon)\subset B_{d_X}(gH_X,R_3) \subset B_{d_\mathfrak m}^g(gH_X,\tilde s/2)$ and $B_{d_\mathfrak m}^g(gH_X,\tilde s/2)\cap A_3^g=\emptyset$. Then, if $O$ is a $\tau$-connected component of $B_{d_X}(H_X,\varepsilon)$, then
\begin{equation}
\label{dentro ss2}
O \subset B^g_{d_\mathfrak m}(gH_X,\tilde s/2) \text{ or else }O \cap B^g_{d_\mathfrak m}[gH_X, s] = \emptyset.
\end{equation} 

$B_{d_X}(H_X,\varepsilon)$ is a $\tau$-open subset of $G/H_X$ and we will show that it is contained in a countable union of disjoint compact sets. 
Let $\{O_i\}_{i\in \Gamma}$ be the $\tau$-connected components of $B_{d_X}(H_X,\varepsilon)$, where $\Gamma$ is a countable set of indexes. 
The idea is to join the connected components in equivalence classes that are in the same $B_{d_\mathfrak m}^g(gH_X,\tilde s/2)$. 
We say that $O_i\sim O_j$ if for every $(x,y)\in O_i\times O_j$, there exist a $\tau$-path $\eta:[a,b]\rightarrow G/H_X$ connecting $x$ and $y$ such that $d_X(\gamma(t_1),\gamma(t_2))\leq R_2$ for every $t_1,t_2\in [a,b]$. 
This implies that if $gH_X\in O_i$, then $O_j\subset B_{d_\mathfrak m}^g(gH_X,s/2)$, that can be refined to $O_j\subset B_{d_\mathfrak m}^g(gH_X,\tilde s/2)$ due to (\ref{bi opcoes}) and (\ref{dentro ss2}).
Observe that $\sim$ is an equivalence relation. 
In fact, reflexivity and symmetry of $\sim$ are immediate. 
In order to see the transitivity of $\sim$, observe that $O_i\sim O_j$ and $O_j\sim O_k$ implies that for every $(x,z)\in O_i\times O_k$, there exist a $\tau$-path connecting $x$ and $z$ in such a way that $d_X(\gamma(t_1),\gamma(t_2))\leq 2R_2<R_1$ for every $t_1,t_2\in [a,b]$. 
Therefore $O_k\subset B_{d_\mathfrak m}(gH_X,s/2)$ for some (any) $gH_X\in O_i$, what implies that $O_k\subset B_{d_\mathfrak m}(gH_X,\tilde s/2)$ due to (\ref{bi opcoes}) and (\ref{dentro ss2}) and we have  that $d_X(\gamma(t_1),\gamma(t_2))\leq R_2$ for every $t_1,t_2\in [a,b]$. 
Therefore $O_i\sim O_k$ and $\sim$ is an equivalence relation.

Observe that if $O_i\cap B_{d_\mathfrak m}^g[gH_X,s]\neq \emptyset$ for some $gH_X\in B_{d_X}(H_X,\varepsilon)$, then $O_i$ must be contained in $B_{d_\mathfrak m}^g(gH_X,\tilde s/2)$ due to (\ref{bi opcoes}) and (\ref{dentro ss2}). 
In particular, if we also have $O_j \cap B_{d_\mathfrak m}^g[gH_X,s]\neq \emptyset$, then $O_i\sim O_j$. 
Therefore, if $O_i\not \sim O_j$ with $O_i\subset B_{d_\mathfrak m}^{g_i}[g_iH_X,s]$ and $O_j\subset B_{d_\mathfrak m}^{g_j}[g_jH_X,s]$, then 
\begin{equation}
\label{compacto nao intercepta}
B_2^{g_i}\cap B_2^{g_j}=\emptyset
\end{equation}
due to (\ref{bi opcoes}) and (\ref{dentro ss2}).

Finally we denote the union $\cup_{i\in \Gamma}O_i$ by $\cup_{l\in \Lambda} U_l$, where each $U_l$ is the union of $\tau$-components of $B_{d_X}(H_X,\varepsilon)$ that are in the same equivalence class defined by $\sim$.
For every $l\in \Lambda$, choose $g_lH_X\in U_l$.  Define $K_l:=B_2^{g_l}$. 
Then $U_l \subset B_{d_\mathfrak m}^{g_l}(g_lH_X,\tilde s/2) \subset K_l$ and $K_l\cap K_m=\emptyset$ whenever $l\neq m$ due to (\ref{compacto nao intercepta}).$\blacksquare$

\section{Paths in $G/H_X$}
\label{caminhos}

In this section we study properties like rectifiability and speed of paths in $(G/H_X,d_X)$ and $(G/H_X,\hat d_X)$. 
We also prove that in some cases, a $d_X$-path in $G/H_X$ is always a $\tau$-path (see Theorem \ref{caminhos coincidem}). 
This theorem is essential in order to prove item (2) of Theorem \ref{lipschitzlocalequivalente}. 

\begin{proposition}
\label{locallispschitzandrectifiable} If $\gamma :[a,b]\rightarrow (M,d)$ is a rectifiable path on a metric space $(M,d)$ and $\rho$ is locally Lipschitz equivalent to $d$, then $\gamma:[a,b]\rightarrow (M,\rho)$ is also rectifiable. Moreover if $\gamma:[a,b]\rightarrow (M,d)$ is Lipschitz, then $\gamma:[a,b]\rightarrow (M,\rho)$ is also Lipschitz.
\end{proposition}

{\it Proof}

\

Suppose that $\gamma:[a,b]\rightarrow (M,d)$ is rectifiable. Notice that $\gamma([a,b])$ is compact. Then $d$ and $\rho$ are Lipschitz equivalent on $\gamma([a,b])$ (see Proposition \ref{compactolipschitzequivalente}). 
Thus there exist a constant $C>0$ such that $\rho(x,y)\leq C.d(x,y)$ for every $x,y\in \gamma([a,b])$ and it follows that $\ell_\rho(\gamma)\leq C\ell_d(\gamma)$. 

The proof of the Lipschitz case is analogous$\blacksquare$

\begin{proposition}
\label{suaveretificavel}
Let $G$ be a Lie group, $(M,d)$ be a differentiable manifold endowed with a metric $d\in \mathcal L(M)$ and $X\subset M$ be a compact subset.
Let $\varphi:G\times M\rightarrow M$ be a smooth action of $G$ on $M$. 
Then every path $\eta:[a,b]\rightarrow (G/H_X,d_X)$ which is continuously differentiable by parts is Lipschitz. 
In particular, $\eta$ is rectifiable. 
These results hold for $\eta:[a,b]\rightarrow (G/H_X,\hat d_X)$ as well.
\end{proposition}

{\it Proof}

\

Due to Theorem \ref{lipschitzimplicalipschitz} and Proposition \ref{locallispschitzandrectifiable}, it is enough to prove the result when $(M,F)$ is a Finsler manifold. 
It is also enough to prove the result only for continuously differentiable curves $\eta$.

Let $\eta:[a,b]\rightarrow G/H_X$ be a continuously differentiable curve.
Then there exist a curve $\tilde \eta:[a,b]\rightarrow G$ which is continuously differentiable by parts and such that $\eta=\pi\circ \tilde \eta$, where $\pi:G\rightarrow G/H_X$ is the natural projection.
In fact, $\pi$ is a submersion.
Then, for every $t\in (a,b)$ and every $p\in \pi^{-1}(\eta(t))$, it is not difficult to find an $\varepsilon>0$ and a continuously differentiable path $\bar \eta:[t-\varepsilon,t+\varepsilon]\rightarrow G$ such that $\pi\circ \bar \eta=\eta\vert_{[t-\varepsilon,t+\varepsilon]}$ and $\bar \eta(t)=p$.
Moreover, if $h\in H_X$, then $\tilde \eta(\cdot):=\bar \eta(\cdot) .h$ satisfies $\pi \circ \tilde \eta=\eta\vert_{[t-\varepsilon,t+\varepsilon]}$ and we can make a ``vertical displacement'' of $\bar \eta$ along the fibers of $\pi$.
For $t\in \{a,b\}$ an analogous continuously differentiable lifting $\bar \eta$ can be made. 
Then we can use this construction and the compactness of $[a,b]$ in order to lift $\eta$ locally and by parts and we can find a curve $\tilde\eta:[a,b]\rightarrow G$ which is continuously differentiable by parts and satisfies $\pi\circ \tilde \eta=\eta$ (We could smooth the vertices of $\tilde \eta$ and the lifting can be made continuously differentiable, but we do not need this fact here).

Then

\[
d_X(\eta(a),\eta(b))\leq
\ell_{d_X}(\eta)=\sup_{\mathcal P}\sum_{i=1}^{n_\mathcal P}d_X(\tilde \eta(t_i)H_X,\tilde \eta(t_{i-1})H_X)
\]
\begin{equation}
\label{somatoriacomprimentocurva1}
=\sup_{\mathcal P}\sum_{i=1}^{n_\mathcal P}d_H(\tilde \eta(t_i)X,\tilde \eta(t_{i-1})X)=C,
\end{equation}
where $\mathcal P$ is the partition $\{a=t_0<t_1<\ldots < t_{n_{\mathcal P}}= b\}$. Without loss of generality, we can suppose that points of $\tilde \eta$ that are not differentiable are in $\mathcal P$. Therefore 
\[
C\leq \sup_{\mathcal P}\sum_{i=1}^{n_\mathcal P}\sup_{z\in X}d_M(\tilde\eta(t_i)z,\tilde \eta(t_{i-1})z)
\leq \sup_{\mathcal P}\sum_{i=1}^{n_\mathcal P}\sup_{z\in X}\int_{t_{i-1}}^{t_i}F\left( \frac{d}{dt} (\tilde \eta(t)z) \right) dt
\]
\begin{equation}
\label{limitantesuperiorcomprimento}
\leq \sup_{z\in X}\int_a^b F\left( \frac{d}{dt} (\tilde \eta(t)z) \right) dt \leq \max_{\substack{z\in X\\ t\in [a,b]}} F\left( \frac{d}{dt} (\tilde \eta(t)z) \right) (b-a)
\end{equation}
and $\eta$ is Lipschitz.

The last statement of the theorem holds because $\ell_{\hat d_X}(\eta)=\ell_{d_X}(\eta)$ for every rectifiable path $\eta:[a,b]\rightarrow (G/H_X,d_X)$ (see (\ref{comprimentoigual})).$\blacksquare$

\

\begin{remark}
\label{killing}
Let $\varphi:G\times M\rightarrow M$ be a smooth action by isometries of a Lie group $G$ on a Finsler manifold $(M,F)$.
Let $v\in \mathfrak g$ and consider the vector field $K_v(p):=\left.\frac{d}{dt}\right\vert_{t=0} \exp (tv).p$ on $M$. 
Then  $K_v$ is a Killing field on $M$, that is, their flow are isometries (see Corollary \ref{nao interessa conceito de isometria 2}).

Notice that 
\[
\left.\frac{d}{dt}\right\vert_{t=t_0} \tilde\eta(t)z=
\left.\frac{d}{dt}\right\vert_{t=t_0} \tilde\eta(t)\tilde\eta^{-1}(t_0)\tilde\eta(t_0)z=K_{\left(\frac{d}{dt}\vert_{t_0}\tilde\eta(t)\tilde\eta^{-1}(t_0)\right)}(\tilde \eta(t_0)z).
\]
Therefore, when $\varphi$ is an action by isometries, the Lipschitz constant of (\ref{limitantesuperiorcomprimento}) is the maximum of the norms of a family of Killing fields on $M$.
\end{remark}  

\begin{remark}
\label{propkillingemgrupos}
Let $\varphi:G\times M\rightarrow M$ be a smooth action by isometries of a Lie group $G$ on $M=G$ endowed with a (left invariant) Riemannian metric $\left<\cdot, \cdot \right>_M$ and $v\in \mathfrak g$.
Then the Killing field $K_v$ can be written as
\begin{equation}
\label{killingemgrupos}
K_v(g)=\left. \frac{d}{dt} \right\vert_{t=0} \exp(tv)g=d(R_g)_e(v).
\end{equation}
\end{remark}

\begin{proposition}
\label{comprimentointegralkilling}
Let $\varphi:G\times M\rightarrow M$ be a continuous action by isometries of a Lie group $G$ on a metric space $(M,d)$. Let $X\subset M$ be a compact subset. 
Let $\eta: \mathbb R\rightarrow G/H_X$ given by $\eta(t)=\exp(tv)H_X$, with $v\in \mathfrak g$. 
Then $\ell_{d_X} (\eta\vert_{[a,b]})=\ell_{d_X}(\eta\vert_{[0,1]})(b-a)$. 
\end{proposition}

\

{\it Proof}

\

First of all, notice that
\[
\ell(\eta_{[a,b]})=\sup_{\mathcal P}\sum_{n=1}^{n_{\mathcal P}} d_H(\exp(t_iv)X,\exp(t_{i-1}v)X)
\]
\[
=\sup_{\mathcal P}\sum_{n=1}^{n_{\mathcal P}} d_H(\exp(cv)\exp(t_iv)X,\exp(cv)\exp(t_{i-1}v)X)=\ell(\eta_{[a+c,b+c]})
\]
for every $c\in \mathbb R$. Therefore we have that $\ell(\eta_{[0,1]})=\ell(\eta_{[0,1/2]})+\ell(\eta_{[1/2,1]})=2\ell(\eta_{[0,1/2]})$. More in general, it is immediate to see that $\ell(\eta_{[0,1]})=2^k\ell(\eta_{[0,2^{-k}]})$. These facts allow us to calculate $\ell(\eta\vert_{[a,b]})$ for any $a<b$. In fact, every closed interval can be obtained as a countable union of closed intervals with measure $1/2^k$, where the intersection between two of these intervals is a point (eventually we need an additional point). For instance
\[
\left[0,\frac{2}{3}\right]=\left\{\left[0,\frac{1}{2}\right]\cup \left[\frac{1}{2},\frac{1}{2}+\frac{1}{2^3}\right]\cup \left[\frac{1}{2}+\frac{1}{2^3},\frac{1}{2}+\frac{1}{2^3}+\frac{1}{2^5}\right]\ldots \right\} \cup \left\{\frac{2}{3}\right\}
\]
implies that 
\[
\ell ( \eta\vert_{\left[ 0,\frac{2}{3} \right]} )
= \ell(\eta\vert_{[0,\frac{1}{2}]}) + \ell(\eta\vert_{[\frac{1}{2},\frac{5}{8}]}) + \ldots 
\]
\[
= \frac{1}{2} \ell(\eta\vert_{[0,1]})+ \frac{1}{8}\ell(\eta\vert_{[0,1]}) + \ldots 
= \frac{2}{3}\ell(\eta\vert_{[0,1]}).
\]

We can generalize the case $[0,2/3]$ easily and we have that
\[
\label{eqcomprimento}
\ell(\eta_{[a,b]})=\ell(\eta_{[0,1]}).(b-a).\blacksquare
\]

\begin{example}
In Proposition \ref{comprimentointegralkilling}, if $(M,d)$ is a differentiable manifold with $d\in \mathcal L (M)$ and $\varphi:G\times M\rightarrow M$ is a smooth action by isometries, then $\ell_{d_X}(\eta_{[0,1]})<\infty$ due to Proposition \ref{suaveretificavel}. 
This is not always the case if $d\not\in \mathcal L(M)$. For instance, let $G$ be the Heisenberg group and consider the action $\varphi:G\times G\rightarrow G$ given by the product of $G$. Let $X=\{e\}\subset G$ be the compact subset and $(V_1,V_2,V_3)$ be a left invariant moving frame of $G$ such that $[V_1,V_2]=V_3$ and $[V_i,V_j]=0$ otherwise. Consider the Carnot-Carathéodory metric on $G$ that are generated by an invariant inner product on the distribution generated by $\{V_1,V_2\}$ (see \cite{Mitch}, \cite{Montg}). Then $H_X=\{e\}$, $(G/H_X,d_X)$ is $G$ with the original Carnot-Carathéodory metric and $t\mapsto \exp(tV_3(e))\{e\}$ is not rectifiable even if it is restricted to an arbitrarily small interval $[a,b]$.
\end{example}

Now we study the speed of a differentiable curve in $(G/H_X,d_X)$. First of all we remember the  following classical result.

\begin{proposition}
\label{retificavelepsilondelta}Let $\eta:[a,b]\rightarrow (M,d)$ be a rectifiable curve on a metric space $(M,d)$. Then for every $\varepsilon>0$ there exist a $\delta >0$ such that if we have a partition $\mathcal P=\{t_0=a<t_1<\ldots<t_{n_{\mathcal P}}=b\}$ satisfying $\vert \mathcal P\vert<\delta$, then 
\[
\ell(\eta)-\sum_{i=1}^{n_{\mathcal P}} d(\eta(t_i),\eta(t_{i-1}))<\varepsilon.
\]
\end{proposition}

\begin{theorem}
\label{integralkillingdiferenciavel}
Let $\varphi:G\times M\rightarrow M$ be a smooth action by isometries of a Lie group $G$ on a differentiable manifold $(M,d)$, $d\in \mathcal L (M)$. 
Let $X\subset M$ be a compact subset. 
If $\eta:\mathbb R\rightarrow G/H_X$ is given by $\eta(t)=\exp(tv)H_X$, then
\[
\lim_{t\rightarrow 0}\frac{d_X(\exp(tv)H_X,H_X)}{\vert t \vert}
=\lim_{t\rightarrow 0}\frac{d_H(\exp(tv)X,X)}{\vert t\vert}
\]
\begin{equation}
\label{formuladerivadadx}
=\lim_{t\rightarrow 0}\frac{\hat d_X(\exp(tv)H_X,H_X)}{\vert t\vert}
=\ell_{d_X}(\eta_{[0,1]})
=\ell_{\hat d_X}(\eta_{[0,1]})
<\infty.
\end{equation}

\end{theorem}

{\it Proof}

\

We have that $\ell_{d_X}(\eta_{[0,1]})<\infty$ due to Proposition \ref{suaveretificavel}. 
Then $\eta\vert_{[0,1]}$ is rectifiable in $(G/H_X,d_X)$, what implies that $\ell_{d_X}(\eta\vert_{[0,1]}) = \ell_{\hat d_X}(\eta\vert_{[0,1]})$.

The first equality of (\ref{formuladerivadadx}) is the definition of $d_X$. 
We claim that it is enough to prove that
\begin{equation}
\label{equacao fundamental velocidade}
\lim_{t\rightarrow 0^+}\frac{d_X(\exp(tv)H_X,H_X)}{t}=\ell_{d_X}(\eta_{[0,1]})
\end{equation}
in order to prove the second and third equalities of (\ref{formuladerivadadx}). In fact, we know that 
\[
\lim_{t\rightarrow 0^+}\frac{d_X(\exp(tv)H_X,H_X)}{t}
=\lim_{t\rightarrow 0^+}\frac{d_X(H_X,\exp(-tv)H_X)}{t}
\]
\begin{equation}
\label{derivadas iguais}
=\lim_{t\rightarrow 0^-}\frac{d_X(H_X,\exp(tv)H_X)}{-t}
=\lim_{t\rightarrow 0}\frac{d_X(\exp(tv)H_X,H_X)}{\vert t\vert},
\end{equation}
where the second equality is due to the $G$-invariance of $d_X$. In addition we have that
\[
d_X(\exp(tv)H_X,H_X)
\leq \hat d_X(\exp(tv)H_X,H_X)
\leq \ell_{d_X}(\eta_{[0,1]})t
\]
for $t>0$, where the second inequality is due to Proposition \ref{comprimentointegralkilling}. But the $G$-invariance of $d_X$ and $\hat d_X$ implies that
\[
d_X(H_X,\exp(-tv) H_X)
\leq \hat d_X(H_X,\exp(-tv) H_X)
\leq \ell_{d_X}(\eta_{[0,1]})t
\]
for every $t>0$, what gives
\begin{equation}
\label{ineq dx dxchapeu comprimento 2}
d_X(\exp(tv)H_X,H_X)
\leq \hat d_X(\exp(tv)H_X,H_X)
\leq \ell_{d_X}(\eta_{[0,1]})\vert t \vert
\end{equation}
for every $t \in \mathbb R$. 
Therefore (\ref{equacao fundamental velocidade}), (\ref{derivadas iguais}) and (\ref{ineq dx dxchapeu comprimento 2}) settle the second and third equations of (\ref{formuladerivadadx}). 

Let us prove (\ref{equacao fundamental velocidade}).

Consider $\eta\vert_{[0,1]}$. Let $\mathcal P=\{0=t_0<t_1<\ldots<t_{n_{\mathcal P}}=1\}$ be a partition of $[0,1]$ and set
\[
\Sigma(\mathcal P)=\sum_{i=1}^{n_{\mathcal P}}d_H(\exp(t_iv)X,\exp(t_{i-1}v)X).
\]

Fix $\varepsilon>0$. Due to Proposition \ref{retificavelepsilondelta}, we can find $\delta \in (0,1/2)$ such that if $\vert \bar{\mathcal P}\vert<\delta$, then $\ell(\eta)-\Sigma(\bar{\mathcal P})<\varepsilon/2$. 
Let $\mu<\delta$ and $\mathcal P_\mu$ be a partition of $[0,1]$ given by $\{ 0=t_0<\mu < 2\mu < \ldots < N.\mu\leq 1 \}$, where $N$ is chosen such that $t^\prime:=1-N.\mu\in [0,\mu)$. Observe that either $t_N=N.\mu=1$ or $t_N=N\mu<t_{N+1}=1$. Then
\[
\frac{\varepsilon}{2}>\ell_{d_X}(\eta\vert_{[0,1]})-\Sigma({\mathcal P}_\mu)=\ell_{d_X}(\eta\vert_{[0,1]}) - \sum_{i=1}^{n_{\mathcal P_\mu}}d_H(\exp(t_iv)X,\exp(t_{i-1}v)X)
\]
\[
=\ell_{d_X}(\eta_{[0,1]})(N\mu+t^\prime)-\sum_{i=1}^N d_H(\exp(i\mu v)X,\exp((i-1)\mu v) X)
\]
\[
-d_H(\exp(v)X,\exp(N\mu v)X)=C.
\]
Notice that $d_H$ is $G$-invariant. Then
\[
C=\ell_{d_X}(\eta\vert_{[0,1]})N\mu-Nd_H(\exp(\mu v)X,X)+\ell_{d_X}(\eta\vert_{[0,t^\prime]})-d_H(\exp(t^\prime v)X,X)
\]
\begin{equation}
\label{todosiguais}
\geq \ell_{d_X}(\eta\vert_{[0,1]})N\mu-Nd_H(\exp(\mu v)X,X),
\end{equation}
where the equality $\ell_{d_X}(\eta\vert_{[0,1]})t^\prime=\ell_{d_X}(\eta\vert_{[0,t^\prime]})$ is due to Proposition \ref{comprimentointegralkilling}. If we divide both sides of (\ref{todosiguais}) by $N\mu$, we get
\[
\left\vert \ell_{d_X}(\eta\vert_{[0,1]})-\frac{d_H(\exp(\mu v)X,X)}{\mu} \right\vert =\ell_{d_X}(\eta\vert_{[0,1]})-\frac{d_H(\exp(\mu v)X,X)}{\mu}
\]
\[
<\frac{\varepsilon}{2N\mu}<\frac{\varepsilon}{2(1-\delta)}<\varepsilon
\]
for every $\mu\in (0,\delta)$. This settles $(\ref{equacao fundamental velocidade})$ and the theorem.$\blacksquare$

\

In Theorem \ref{funcaoF}, the ($G$-invariant) Finsler metric $F$ correspondent to $(G/H_X,\hat d_X)$ is obtained from
\[
F(H_X,\bar v)=\lim_{t\rightarrow 0} \frac{d_X(c(t),H_X)}{\vert t \vert},
\]
where $c:(-\varepsilon,\varepsilon) \rightarrow G/H_X$ is an arbitrary curve such that $c(0)=H_X$ and $c^\prime(0)=\bar v$. 
This is like a reciprocal of Theorem \ref{velocidade Finsler}.
But we need that $F(H_X,\bar v)$ does not depend on the choice of $c$. The next lemma does the job.

\begin{lemma}
\label{mesmavelocidade}
Let $M$ be a differentiable manifold and denote its topology by $\tau_M$. 
Let $d$ be a metric on $M$ such that for every $p\in M$, there exist a $\tau_M$-neighborhood $V$ of $p$ such that $d\vert_{V\times V}$ is Lipschitz equivalent to a Finsler metric on $V$ ($d$ does not need to induce $\tau_M$). 
Consider $p\in M$, $v\in T_pM$ and suppose that there exist a curve $c:(-\varepsilon,\varepsilon)\rightarrow M$ such that $c(0)=p$, $c^\prime(0)=v$ and in addition such that the speed
\[
\lim_{t\rightarrow 0}\frac{d(c(t),p)}{\vert t\vert}
\]
exists. If $\eta:(-\varepsilon,\varepsilon)\rightarrow M$ is another curve such that $\eta(0)=p$ and $\eta^\prime(0)=v$, then
\[
\lim_{t\rightarrow 0}\frac{d(\eta(t),p)}{\vert t\vert} =
\lim_{t\rightarrow 0}\frac{d(c(t),p)}{\vert t\vert }.
\]
\end{lemma}

{\it Proof}

\

Endow a sufficiently small $\tau_M$-neighborhood $V$ of $p$ with a Euclidean metric $\left<\cdot,\cdot\right>$ and let $(V,\psi=(x_1,\ldots,x_n))$ be a coordinate system of $V$ such that $\psi(p)=(0,\ldots,0)$ and  $\left<\partial/\partial x_i,\partial/\partial x_j\right>=\delta_{ij}$. Without loss of generality, we can also suppose that $V$ is $\left< \cdot, \cdot \right>$-convex and that the distance function $d_{\text{flat}}$ induced by $\left<\cdot,\cdot \right>$ is Lipschitz equivalent to $d$ on $V$.

We can restrict $c$ and $\eta$ in such a way that their images are contained in $V$. Write $c(t)=(c_1(t),\dots,c_n(t))$ and $v=(v_1,\ldots,v_n)$ with respect to the coordinate system $\psi$. 
If we represent $c$ by its Taylor polynomial and its remainder term, then $c(t)=(v_1.t+f_1(t), v_2.t+f_2(t),\ldots v_n.t+f_n(t))$, where $f_i(t)=O(t^2)$ for every $i$.
Likewise the Taylor polynomial of $\eta$  is given by $\eta(t)=(v_1.t+u_1(t),v_2.t+u_2(t),\ldots v_n.t+u_n(t))$, where $u_i(t)=O(t^2)$ for every $i$. Then we have the estimates $d_{\text{flat}}( c(t),\eta(t))=O(t^2)$ and $d(\eta(t),c(t))=O(t^2)$ because $d_{\text{flat}}$ and $d$ are Lipschitz equivalent. 
Due to the triangle inequality, we have that $\vert d(p,c(t))-d(p,\eta(t))\vert =O(t^2)$ and the result follows.$\blacksquare$

\

Proposition \ref{dxcontinua} implies that $\tau$-paths in $G/H_X$ are $d_X$-paths. We conclude this section presenting a type of converse.

\begin{theorem}
\label{caminhos coincidem}
Let $G$ be a Lie group, $M$ be a differentiable manifold endowed with a metric $d$, $X$ be a compact subset of $M$ and $\varphi:G\times M\rightarrow M$ be a smooth action by isometries. Let $I\subset \mathbb R$ be an interval and $\gamma:I\rightarrow G/H_X$ be a $d_X$-path. Then $\gamma:I\rightarrow G/H_X$ is a $\tau$-path.
\end{theorem}

This theorem is important to prove the second item of Theorem \ref{lipschitzlocalequivalente} and it uses the study of the geometry of $(G/H_X,d_X)$ which we made in Theorem \ref{separacao bolas}.

Before the proof of Theorem \ref{caminhos coincidem} we remember some results and prove some lemmas.

\begin{theorem}[Sierpi\'nski theorem]
\label{sierpinski}
Let $X$ be a continuum (compact, connected and Hausdorff topological space). If $\{C_i\}_{i\in \mathbb N}$ is a pairwise disjoint countable covering of $X$ by closed subsets, then $C_i=X$ for some $i\in \mathbb N$. 
\end{theorem}

{\it Proof}

\

See \cite{Sierpinski}.

\

\begin{proposition}
\label{fraco e compacto forte}
Let $M$ be a non-empty set and consider two Hausdorff topologies $\tau_c$ and $\tau_f$ on $M$ such that $\tau_c\subset \tau_f$. 
Suppose that $x_i \stackrel{\tau_c}{\rightarrow} x$ and that $\cup_{i\in \mathbb N} \{x_i\}$ is contained in a $\tau_f$-sequentially compact subset of $M$. Then $x_i\stackrel{\tau_f}{\rightarrow} x$.
\end{proposition}

{\it Proof:} Immediate.

\begin{corollary}
\label{contido em compacto continua}
Suppose we are in the hypotheses of Theorem \ref{caminhos coincidem} and that a $d_X$-path $\gamma:I\rightarrow M$ is contained in a $\tau$-compact subset. Then $\gamma$ is $\tau$-continuous.
\end{corollary}

{\it Proof}

\

We will prove that if $x_i\rightarrow x$ in $I$, then $\gamma(x_i)\stackrel{\tau}{\rightarrow} \gamma(x)$ in $M$. 
We have that $\{\gamma(x_i)\}$ is contained in a $\tau$-compact subset and $\gamma(x_i) \stackrel{d_X}{\rightarrow} \gamma (x)$. Therefore $\gamma (x_i)\stackrel{\tau}{\rightarrow} \gamma (x)$ due to Proposition \ref{fraco e compacto forte} and $\gamma$ is $\tau$-continuous. $\blacksquare$ 




\

{\it Proof of Theorem \ref{caminhos coincidem}}

\

For the sake of simplicity, suppose that $I=[a,b]$ (The proof is local, then it works for all types of intervals). 
Fix $t\in (a,b)$ (The case $t\in \{a,b\}$ will be seen afterwards). 
From Theorem \ref{separacao bolas}, there exist an $\varepsilon >0$ such that $B_{d_X}(\gamma(t),\varepsilon)$ is contained in a countable union of pairwise disjoint $\tau$-compact subsets $\cup_{i\in \Lambda} K_i$.
Due to the $d_X$-continuity of $\gamma$, there exist a $\delta>0$ such that $[t-\delta,t+\delta]\subset [a,b]$ and $\gamma([t-\delta,t+\delta])\subset B_{d_X}(\gamma(t),\varepsilon)$.
We will prove that $\eta:=\gamma\vert_{[t-\delta,t+ \delta]}$ is $\tau$-continuous.

We claim that $\eta^{-1}(K_j)$ is closed for every $j\in \Lambda$.
We will prove that $\eta^{-1}(K_j)$ contains its accumulation points.
Suppose that $x$ is an accumulation point of $\eta^{-1}(K_j)$ and consider a sequence $x_i \rightarrow x$ in $\eta^{-1}(K_j)$.  
Then $\eta(x_i)\stackrel{d_X}{\rightarrow} \eta(x)$ and $\{\eta(x_i)\}\subset K_j$ holds, what implies that $\eta (x_i)\stackrel{\tau}{\rightarrow} \eta(x)$ due to Proposition \ref{fraco e compacto forte}. Thus $\eta(x)\in K_j$ and $x\in \eta^{-1}(K_j)$.

Notice that $\{\eta^{-1}(K_j)\}_{j\in \Lambda}$ is a cover of $[t-\delta,t+\delta]$ by pairwise disjoint closed subsets.
Then only one $\eta^{-1}(K_j)$ is non-empty (see Theorem \ref{sierpinski}) and $\eta([t-\delta,t+\delta])\subset K_j$. 
Then $\eta\vert_{[t-\delta,t+\delta]}$ is $d_X$-continuous and contained in a $\tau$-compact subset what implies that $\eta\vert_{[t-\delta,t+\delta]}$ is $\tau$-continuous due to Corollary \ref{contido em compacto continua}.

The proof for the case $t\in \{a,b\}$ is analogous. 

This proof holds for any type of intervals, what settles the theorem.$\blacksquare$

\section{Intrinsic induced Hausdorff metrics are Finsler}
\label{secaofinsler}

In this section we prove the following theorem: Let $G$ be a Lie group, $M$ be a differentiable manifold endowed with a metric $d\in \mathcal L(M)$ and $X\subset M$ be a compact subset. 
Let  $\varphi:G\times M\rightarrow M$ be a smooth action by isometries of $G$ on $M$. 
Then $(G/H_X,\hat d_X)$ is a Finsler manifold (see Theorem \ref{intrinsecaefinsler}). 
In particular, a homogeneous space endowed with an intrinsic $G$-invariant metric $d\in \mathcal L(M)$ is a Finsler manifold (see Corollary \ref{intrinsecafinsler2}).

The next theorem is essential in order to prove Theorem \ref{intrinsecaefinsler}.

\begin{theorem}
\label{lipschitzlocalequivalente}
Let $\varphi:G\times M\rightarrow M$ be a smooth action by isometries of the Lie group $G$ on a differentiable manifold $(M,d\in \mathcal L(M))$ and let $X\subset M$ be a compact subset. 
Consider an arbitrary Finsler metric $F$ on $G/H_X$ and denote its distance function by $d_F$. 
\begin{enumerate}

\item Consider $gH_X\in G/H_X$. 
Then there exist a $\tau$-neighborhood $O$ of $gH_X$ such that $d_F\vert_{O\times O}$, $d_X\vert_{O\times O}$ and $\hat d_X\vert_{O\times O}$ are Lipschitz equivalent.

\item If $\tau_X$ is the topology induced by $d_X$ and $\hat \tau_X$ is the topology induced by $\hat d_X$, then $\tau_X\subset \hat \tau_X=\tau$.

\end{enumerate}

\end{theorem}

\

Before proving this theorem, we prove some preliminary results.

\

\begin{lemma}
\label{cobreintervalo}
Consider $(a-\varepsilon,a+\varepsilon)\subset \mathbb R$, with $\varepsilon,a-\varepsilon>0$. Then there exist a $\delta>0$ such that 
\[
\bigcup_{n\in \mathbb N} \left(\frac{a-\varepsilon}{n},\frac{a+\varepsilon}{n}\right)\supset (0,\delta).
\]
\end{lemma}

{\it Proof}

\

It is enough to show that there exist a $N\in \mathbb N$ such that $(a+\varepsilon)/(n+1)>(a-\varepsilon)/n$ for every $n\geq N$. 
But it is straightforward that this is true for $n>(a-\varepsilon)/(2\varepsilon)$.
Therefore if we choose any $N>(a-\varepsilon)/(2\varepsilon)$, then
\[
\bigcup_{n\geq N} \left(\frac{a-\varepsilon}{n},\frac{a+\varepsilon}{n}\right)= \left(0,\frac{a+\varepsilon}{N}\right).\blacksquare
\]

\

\begin{lemma}
\label{distanciaaorigemlipschitz}

\

\begin{itemize}
\item Let $\varphi:G\times M\rightarrow M$ be a continuous action by isometries of a Lie group $G$ on a metric space $M$ and let $X$ be a compact subset of $M$. 
Consider an arbitrary Finsler metric $F$ on $G/H_X$ and denote its distance function by $d_F$. 
Then there exist a constant $C>0$ and a $\tau$-neighborhood $O$ of $H_X\in G/H_X$ such that
\begin{equation} 
\label{dfdx}
d_F(H_X,gH_X)\leq Cd_X(H_X,gH_X)
\end{equation} 
for every $gH_X\in O$.

\item  Let $\varphi:G\times M\rightarrow M$ be a smooth action of a Lie group $G$ on a Finsler manifold $(M,F^\prime)$ and let $X$ be a compact subset of $M$. Then there exist a constant $C^\prime >0$ and a $\tau$-neighborhood $O^\prime$ of $H_X\in G/H_X$ such that 
\begin{equation}
\label{hatdxdf}
\hat d_X(H_X,gH_X)\leq C^\prime d_F(H_X,gH_X)
\end{equation} 
for every $gH_X\in O^\prime$.
\end{itemize}
\end{lemma}

{\it Proof}

\

Consider the general setting as in \ref{mhx}. Fix $r^\prime\in (0,r)$ and \ an $\varepsilon >0$ such \ that \   $[r^\prime-\varepsilon,r^\prime+\varepsilon] \subset (0,r)$. 
Observe that $d_F$ and $d_\mathfrak m$ are Lipschitz equivalent on $B_{d_\mathfrak m}[H_X,r^\prime + \varepsilon]$ (see Proposition \ref{compactolipschitzequivalente}).

Let us prove that there exist $C>0$ and a $\tau$-neighborhood $O$ of $H_X$ such that $d_F(H_X,gH_X)\leq Cd_X(H_X,gH_X)$ for every $gH_X\in O$.

Set
\[
\tilde C=\min_{v\in A_{d_\mathfrak m} [0,r^\prime-\varepsilon,r^\prime+\varepsilon]} d_X(exp(v)H_X,H_X),
\]
where $A_{d_\mathfrak m} [0,r^\prime-\varepsilon,r^\prime+\varepsilon] \subset \mathfrak m$. Then
\[
\frac{d_\mathfrak m(\exp(v)H_X,H_X)}{r}
\leq 1
\leq \frac{d_X(exp(v)H_X,H_X)}{\tilde C}
\]
for every $v\in A_{d_\mathfrak m}(0,r^\prime-\varepsilon,r^\prime+\varepsilon)$, what implies that
\[
d_\mathfrak m(\exp(v)H_X,H_X)\leq C d_X(exp(v)H_X,H_X)
\]
where $C=r/\tilde C$. Observe that if $k\in \mathbb N$, then 
\[
d_X(H_X,\exp(v)H_X )\leq k.d_X(H_X,\exp(v/k)H_X)
\] 
for every $v\in A_{d_\mathfrak m}(0,r^\prime-\varepsilon,r^\prime+\varepsilon)$ due to the $G$-invariance of $d_X$ and the triangle inequality. Moreover 
\[
d_{\mathfrak m}(H_X,\exp(v)H_X)= k.d_{\mathfrak m}(H_X,\exp(v/k)H_X)
\]
holds. 
Therefore we have that 
\begin{equation}
\label{comprimentopositivo}
d_{\mathfrak m}(H_X,\exp (tw)H_X)\leq C d_X(H_X,\exp (tw)H_X)
\end{equation}
for every $v\in S_{d_\mathfrak m}(0,1)$ and every
\[
t \in \bigcup_{k\in \mathbb N} \left(\frac{r^\prime-\varepsilon}{k}, \frac{r^\prime+\varepsilon}{k}\right).
\]
Due to Lemma \ref{cobreintervalo}, there exist a $\delta>0$ such that 
\[
d_{\mathfrak m}(H_X,\exp (tw)H_X)\leq Cd_X(H_X,\exp (tw)H_X)
\]
for every $t\in (0,\delta)$ and every $w\in S_{d_\mathfrak m}(0,1)$, that is, for every $tw\in B_{d_\mathfrak m}(0,\delta)$. 
Consequently if we choose $O=B_{d_\mathfrak m}(H_X,\delta)$, then item (1) is settled because $d_F$ and $d_\mathfrak m$ are Lipschitz equivalent on $B_{d_\mathfrak m}(H_X,\delta)$.

Now we prove that $\hat d_X\leq C^\prime d_F$ in a $\tau$-neighborhood of $H_X$. Consider $u\in B_{d_{\mathfrak m}}(0,r/2)-\{0\}$ and $\eta:[0,\Vert u\Vert_{d_{\mathfrak m}}]\rightarrow G/H_X$, $\eta(t)=\exp(tu/\Vert u\Vert_{d_{\mathfrak m}})H_X$. Then
\[
\hat d_X(H_X,\exp(u)H_X)\leq \ell_{\hat d_X}(\eta)\leq 
\max_{\substack{z\in X\\ t\in [0,\Vert u \Vert_{d_{\mathfrak m}}]}} F\left( \frac{d}{dt} (\exp(tu/\Vert u\Vert_{d_{\mathfrak m}})z) \right) \Vert u \Vert_{d_{\mathfrak m}}=C_1,
\]
where the last inequality is calculated as in (\ref{limitantesuperiorcomprimento}). But $\Vert u \Vert_{d_{\mathfrak m}}=d_{\mathfrak m}(H_X,\exp(u)H_X)$. Then
\begin{equation}
\label{formulaintrinsecamenosfina}
C_1\leq  \max_{\substack{z\in X\\ t\in [0,r/2]\\ v\in S_{d_{\mathfrak m}}(0,1)}} F\left( \frac{d}{dt} (\exp(tv)z) \right) d_{\mathfrak m}(H_X,\exp(u)H_X)\leq \tilde C d_{\mathfrak m}(H_X,\exp(u)H_X).
\end{equation}
for every $u\in B_{d_\mathfrak m}(0,r/2)$. Now we use the Lipschitz equivalence of $d_\mathfrak m$ and $d_F$ on $B_{d_\mathfrak m}(H_X,r/2)$ to conclude that there exist a $d_F$-neighborhood $O^\prime=B_{d_\mathfrak m}(H_X,r/2)$ of $H_X \in G/H_X$ and a $C^\prime >0$ such that $\hat d_X(H_X,gH_X)\leq C^\prime d_F(H_X,gH_X)$ for every $gH_X\in O^\prime$.$\blacksquare$

\

\begin{lemma}
\label{topologia local mais fina} Let $M$ be a non-empty set and consider two topologies $\tau_c$ and $\tau_f$ on $M$. 
Suppose that for every $p\in M$ there exist a $\tau_f$-neighborhood $V$ of $p$ such that $\tau_f\vert_V$ is finer than $\tau_c\vert_V$. Then $\tau_f$ is finer than $\tau_c$.
\end{lemma}

{\it Proof}

\

Consider $O_c$ a $\tau_c$-open subset of $M$ and $p\in O_c$. 
We will find a $\tau_f$-open subset $O_f\ni p$ such that $O_f\subset O_c$.
Observe that $O_c\cap V$ is $\tau_c\vert_V$-open.
Then $O_c\cap V$ is also $\tau_f\vert_V$-open, because $\tau_f\vert_V$ is finer than $\tau_c\vert_V$.
Therefore $O_c \cap V$ can be written as $A \cap V$, where $A$ is an $\tau_f$-open subset of $M$, what implies that $O_f:=O_c\cap V\subset O_c$ is $\tau_f$-open.$\blacksquare$

\

The following lemma helps to show that $\hat d_X$ is finer than $\tau$ in $G/H_X$. 
The lemma is applied with $d=d_X$ and $\rho = d_F$.

\begin{lemma}
\label{d chapeu mais fina}
Let $M$ be a non-empty set endowed with metrics $d$ and $\rho$ such that
\begin{itemize}
\item paths in $(M,d)$ are paths in $(M,\rho)$;
\item for every $p\in M$, there exist a $\rho$-neighborhood $V$ of $p$ such that $d\vert_V$ is Lipschitz equivalent to $\rho\vert_V$.
\end{itemize}.  
Then $\hat d$ is finer than $\rho$.
\end{lemma}

{\it Proof}

\

Fix $p \in M$. For an $\varepsilon >0$ we will prove that there exist $\delta >0$ such that $B_{\hat d}(p,\delta)\subset B_\rho(p,\varepsilon)$.
We can consider without loss of generality that $B_\rho[p,\varepsilon]\subset V$. 
Due to the definition of $V$, there exist a $C>0$ such that $\rho(x,y)\leq C d(x,y)$ for every $x,y\in V$. 
If $\eta$ is a $d$-path that does not remain in $V$, then it is also a $\rho$-path that does not remain in $B_\rho(p,\varepsilon)$ and we have that $\ell_d(\eta)\geq \varepsilon /C$.

We claim that $B_{\hat d}(p,\varepsilon /C) \subset B_\rho (p,\varepsilon)$.
If $q \in B_{\hat d}(p,\varepsilon /C)$, then there exist a rectifiable $d$-path $\gamma$ connecting $p$ and $q$ such that $\ell_{\hat d}(\gamma)=\ell_d(\gamma)<\varepsilon/C$ what implies that $\gamma$ does not leave $V$. Then $\ell_\rho(\gamma) < \varepsilon$ and $q \in B_\rho(p,\varepsilon)$, what settles the lemma.$\blacksquare$

\begin{remark}
\label{diferenciavel e continua}

\

\begin{itemize}

\item 
Consider the smooth map
\[
d\phi:T(G\times G/H_X)\rightarrow T(G/H_X).
\]
Represent the smooth curves $\eta:I\rightarrow G\times G/H_X$ by $\eta(t)=(\eta_1(t),\eta_2(t))$ where $\eta_1$ is a curve on $G$ and $\eta_2$ is a curve on $G/H_X$. When we restrict $d\phi$ to the vectors of type $(0,v)$, that is, directional derivatives along curves with $\eta_1$ constant, then 
\begin{equation}
\label{diferenciabilidadedphi}
d\phi(g,hH_X,(0,v)) = (ghH_X,d(\phi_g)_{(hH_X)}(v)).
\end{equation}
In particular, $d(\phi_g)_{hH_X}(v)$ depends continuously on $g$, $hH_X$ and $v$.

\item Suppose that we are in the settings \ref{mhx}. Consider the restriction 
\[
\tilde \phi: \exp B_{d_\mathfrak m}(0,r) \times B_{d_\mathfrak m}(H_X,r) \rightarrow G/H_X
\] 
of $\phi$ and observe that $\exp B_{d_\mathfrak m}(0,r)$ is a submanifold of $G$ which is invariant by inversion. 
We can restrict $d \tilde \phi$ to the points of type $(g^{-1},gH_X,0,v)$ and we get
\[
d\tilde \phi(g^{-1},gH_X,0,v)=(H_X,d(\tilde \phi_{g^{-1}})_{gH_X}(v)).
\]
In particular $\xi:T(B_{d_\mathfrak m}(H_X,r)) \rightarrow T_{H_X}(G/H_X)$ given by 
\begin{equation}
\label{invariante continua}
\xi (gH_X,v) = d(\tilde \phi_{g^{-1}})_{gH_X} (v)=d( \phi _{g^{-1}})_{gH_X} (v),
\end{equation}
is continuous.

\end{itemize}

\end{remark}

\

{\it Proof of Theorem \ref{lipschitzlocalequivalente}}

\

Denote by $\tilde O$ a $\tau$-neighborhood of $H_X\in G/H_X$ with compact closure such that (\ref{dfdx}) and (\ref{hatdxdf}) holds. First of all, we prove the theorem for the case where $(M,d)=(M,F)$ is a Finsler manifold. 

Due to the $G$-invariance of $d_X$ and $\hat d_X$ on $G/H_X$, it is enough to prove that there exist a $\tau$-neighborhood $O$ of $H_X\in G/H_X$ and $C>0$ such that 
\begin{equation}
\label{compara1}
\hat d_X(g_1H_X,g_2H_X)\leq Cd_F(g_1H_X,g_2H_X)
\end{equation}
and 
\begin{equation}
\label{compara2}
d_F(g_1H_X,g_2H_X)\leq Cd_X(g_1H_X,g_2H_X)
\end{equation}
hold for every $g_1H_X,g_2H_X \in O$. Remember that the inequality $d_X\leq \hat d_X$ always holds. 

Let us prove (\ref{compara1}). 
Due to (\ref{diferenciabilidadedphi}), $d(\phi_g)_{hH_X}(v)$ depends continuously on $g$, $hH_X$ and $v$. Then there exist a $(G\times \tau)$-neighborhood $O_1 \times O_2 \subset G\times G/H_X$ of $(e,H_X)$ such that:
\begin{enumerate}
\item $\phi(O_1,O_2)\subset \tilde O$;
\item $\sup\limits_{\Vert v\Vert_{\mathfrak m}=1}\Vert d(\phi_g)_{hH_X}(v)\Vert_{\mathfrak m}\leq 2$ for every $(g,hH_X)\in O_1\times O_2$;
\item $O_1$ is a symmetric neighborhood of $e$;
\item $\phi(O_1,H_X)\supset O_2$;
\item $O_2$ is $d_{\mathfrak m}$-convex.
\end{enumerate}
The properties above can be used as steps to construct $O_1$ and $O_2$ in such a way that at the end of the process all the properties hold.

If $g_1H_X, g_2H_X \in O_2$ and $\gamma:[0,1]\rightarrow (O_2,d_{\mathfrak m})$ is the geodesic connecting $g_1H_X$ and $g_2H_X$ (an Euclidean segment), then $\phi_{g_1^{-1}}\circ \gamma$ is smooth curve in $\tilde O$ (not necessarily an $d_{\mathfrak m}$ geodesic) connecting $H_X$ and $g_1^{-1}g_2 H_X$ such that
\[
d_{\mathfrak m}(H_X,g_1^{-1}g_2H_X)\leq \ell_{d_{\mathfrak m}}(\phi_{g_1^{-1}}\circ \gamma)=\int_0^1 \Vert (\phi_{g_1^{-1}}\circ \gamma)^\prime(t)\Vert_{\mathfrak m}dt
\]
\[
=\int_0^1 \Vert d(\phi_{g_1^{-1}})_{\gamma(t)}(\gamma^\prime(t))\Vert_{\mathfrak m}dt\leq 2\int_0^1 \Vert \gamma^\prime(t)\Vert_{\mathfrak m}dt =2d_{\mathfrak m}(g_1H_X,g_2H_X).
\]
Therefore we have that
\[
\hat d_X(g_1H_X,g_2H_X)=\hat d_X(H_X,g_1^{-1}g_2H_X)\leq C.d_{\mathfrak m}(H_X,g_1^{-1}g_2H_X)
\]
\[
 \leq 2C d_{\mathfrak m}(g_1H_X,g_2H_X),
\]
and (\ref{compara1}) is proved for  $g_1H_X,g_2H_X\in O_2$.

In order to prove (\ref{compara2}), let $O_1^\prime \times O_2^\prime \subset O_1 \times O_2$ be a $(G\times \tau)$-neighborhood of $(e,H_X)$ such that:
\begin{enumerate}
\item $\phi(O_1^\prime,O_2^\prime)\subset O_2$;
\item $O_1^\prime$ is a symmetric neighborhood of $e$;
\item $\phi(O_1^\prime,H_X)\supset O_2^\prime$.
\end{enumerate}
The property $\sup\limits_{\Vert v\Vert_{\mathfrak m}=1}\Vert d(\phi_g)_{hH_X}(v)\Vert_{\mathfrak m}\leq 2$ for every $(g,hH_X)\in O_1^\prime \times O_2^\prime$ is automatically satisfied. In addition, we do not use the condition that $O_2^\prime$ is $d_{\mathfrak m}$-convex.

Let $g_1H_X,g_2H_X\in O_2^\prime$. Then $H_X,g_1^{-1}g_2H_X \in O_2$ due to the properties of $O_1^\prime\times O_2^\prime$. But $O_2$ is $d_{\mathfrak m}$-convex and the and we can consider the $d_{\mathfrak m}$-geodesic $\gamma:[0,1]\rightarrow (O_2,d_{\mathfrak m})$ connecting $H_X$ and $g_1^{-1}g_2H_X$. Then $\phi_{g_1}\circ \gamma$ is smooth curve on $\tilde O$ (not necessarily an $d_{\mathfrak m}$ geodesic) connecting $g_1H_X$ and $g_2 H_X$ such that
\[
d_{\mathfrak m}(g_1 H_X,g_2H_X)\leq \ell_{d_{\mathfrak m}}(\phi_{g_1}\circ \gamma)=\int_0^1 \Vert (\phi_{g_1}\circ \gamma)^\prime(t)\Vert_{\mathfrak m}dt
\]
\[
=\int_0^1 \Vert d(\phi_{g_1})_{\gamma(t)}(\gamma^\prime(t))\Vert_{\mathfrak m}dt\leq 2\int_0^1 \Vert \gamma^\prime(t)\Vert_{\mathfrak m}dt =2d_{\mathfrak m}(H_X,g_1^{-1}g_2H_X).
\]
Therefore 
\[
d_{\mathfrak m}(g_1H_X,g_2H_X)\leq 2d_{\mathfrak m}(H_X,g_1^{-1}g_2H_X)\leq 2Cd_X(H_X,g_1^{-1}g_2H_X)
\]
\[
= 2C d_X(g_1H_X,g_2H_X)
\]
for every $g_1H_X,g_2H_X\in O_2^\prime$. This settles item (1) of the theorem with $O=O_2^\prime$ for the case where $d$ is a Finsler metric because $d_F\vert_{O_2^\prime \times O_2^\prime}$ is Lipschitz equivalent to $d_{\mathfrak m}\vert_{O_2^\prime \times O_2^\prime}$ (see Proposition \ref{compactolipschitzequivalente}). 

In order to prove item (2) of the theorem for the case where $(M,d)$ is Finsler, notice that $\tau_X\subset \tau$ due to Proposition \ref{dxcontinua}. 
The relationship $\hat \tau_X\subset \tau$ is consequence of item (1) and Lemma \ref{topologia local mais fina} (with $\tau_c=\hat \tau_X$ and $\tau_f=\tau$). 
The relationship $\tau \subset \hat \tau_X$ is consequence of item (1), Theorem \ref{caminhos coincidem} and Lemma \ref{d chapeu mais fina} (In Lemma  \ref{d chapeu mais fina}, set $d=d_X$ and $\rho=d_F$).
This settles the case where $(M,d)$ is Finsler.

For the case $(M,d)$ with $d\in \mathcal L(M)$, fix an arbitrary Finsler metric $F^\prime$ on $M$ with distance function $d_{F^\prime}$. 
Then $d_X$ is locally Lipschitz equivalent to $(d_{F^\prime})_X$ (Theorem \ref{lipschitzimplicalipschitz}), what implies that for every $gH_X\in G/H_X$, there exist a $d_X$-neighborhood $O^\prime$ of $gH_X$ and positive constants $c^\prime, C^\prime > 0$ such that 
\begin{equation}
\label{desigualdadelipschitz1}
c^\prime (d_{F^\prime})_X(g_1H_X,g_2H_X)\leq d_X(g_1H_X,g_2H_X) \leq C^\prime (d_{F^\prime})_X(g_1H_X,g_2H_X) 
\end{equation} 
for every $g_1H_X,g_2H_X\in O^\prime$. Now we use the first part of the proof (the case where $(M,d)$ is Finsler) and the fact that $\tau$ is finer than $\tau_X$ in order to find a $\tau$-neighborhood $O\subset O^\prime$ of $gH_X$ and positive constants $c,C>0$ such that 
\begin{equation}
\label{desigualdadelipschitz2}
c d_F(g_1H_X,g_2H_X)\leq (d_{F^\prime})_X(g_1H_X,g_2H_X) \leq C d_F(g_1H_X,g_2H_X)
\end{equation} 
for every $g_1H_X,g_2H_X\in O$. Now we combine (\ref{desigualdadelipschitz1}) and (\ref{desigualdadelipschitz2}) in order to get
\[
cc^\prime d_F(g_1H_X,g_2H_X)\leq d_X(g_1H_X,g_2H_X) \leq CC^\prime d_F(g_1H_X,g_2H_X)
\]
for every $g_1H_X,g_2H_X\in O$. Then $d_X\vert_{O\times O}$ is Lipschitz equivalent to $d_F\vert_{O\times O}$.

The existence of a $d_F$-neighborhood $O$ of an arbitrary $gH_X\in G/H_X$ such that $d_F\vert_{O\times O}$ is Lipschitz equivalent to $\hat d_X\vert_{O\times O}$ is analogous: 
Just replace $d_X$ and $(d_{F^\prime})_X$ of (\ref{desigualdadelipschitz1}) by $\hat d_X$ and $(\hat d_{F^\prime})_X$ respectively and replace $(d_{F^\prime})_X$ of (\ref{desigualdadelipschitz2}) by $(\hat d_{F^\prime})_X$.

Finally observe that the relationships $\tau_X\subset \hat \tau_X=\tau$ follow from Theorem \ref{lipschitzimplicalipschitz}.$\blacksquare$ 

\begin{remark}
In Example \ref{contralocalLipschitz} we have that $\tau_X\underset{\neq}{\subset}\tau$.
\end{remark}

Let $G$ be a Lie group, $(M,d\in \mathcal L(M))$ be a differentiable manifold, $X$ be a compact subset of $M$ and $\varphi: G \times M\rightarrow M$ be a smooth action by isometries of $G$ on $M$. The next step is to construct the function $F: T(G/H_X)\rightarrow \mathbb R$ that will be the Finsler function correspondent to $\hat d_X$ (Theorem \ref{funcaoF}).

\begin{corollary}
\label{naodependedacurva}
Let $\varphi:G\times M\rightarrow M$ be a smooth action by isometries of a Lie group $G$ on a differentiable manifold $(M,d\in \mathcal L(M))$ and let $X\subset M$ be a compact subset. 
Let $\bar v:=v+\mathfrak h_X\in T_{H_X}(G/H_X)$ and $\eta :(-\varepsilon,\varepsilon)\rightarrow G/H_X$ such that $\eta (0)=H_X$ and $\eta^\prime(0)=\bar v$. Then
\[
\lim_{t\rightarrow 0}\frac{d_X(H_X,\eta (t))}{\vert t\vert}=\lim_{t\rightarrow 0}\frac{d_X(H_X,\exp(tv)H_X)}{\vert t\vert}
\] 
and
\[
\lim_{t\rightarrow 0}\frac{\hat d_X(H_X,\eta (t))}{\vert t\vert}=\lim_{t\rightarrow 0}\frac{\hat d_X(H_X,\exp(tv)H_X)}{\vert t\vert}.
\] 
In particular
\[
\lim_{t\rightarrow 0}\frac{d_X(H_X,\exp(tw)H_X)}{\vert t\vert}=\lim_{t\rightarrow 0}\frac{d_X(H_X,\exp(tv)H_X)}{\vert t\vert}
\] 
and
\[
\lim_{t\rightarrow 0}\frac{\hat d_X(H_X,\exp(tw)H_X)}{\vert t\vert}=\lim_{t\rightarrow 0}\frac{\hat d_X(H_X,\exp(tv)H_X)}{\vert t\vert}
\] 
whenever $v-w\in \mathfrak h_X$.
\end{corollary}

\

{\it Proof}

\

It is a direct consequence of Theorem \ref{lipschitzlocalequivalente} and Lemma \ref{mesmavelocidade} (with $c(t)=\exp(tv)H_X$).$\blacksquare$

\begin{theorem}
\label{funcaoF} 
Let $\varphi:G\times M\rightarrow M$ be a smooth action by isometries of a Lie group $G$ on a differentiable manifold $(M,d\in \mathcal L (M))$ and let $X$ be a compact subset of $M$.
Let $\phi: G \times G/H_X \rightarrow G/H_X$ be the natural action of $G$ on $G/H_X$.
Define the function $F_{H_X}:T_{H_X}(G/H_X)\rightarrow \mathbb R$ by
\[
F_{H_X}(\bar v)=\lim_{t\rightarrow 0}\frac{d_X(c(t),H_X)}{\vert t \vert},
\] 
where $c:(-\varepsilon,\varepsilon)\rightarrow G/H_X$ is any curve such that $c(0)=H_X$ and $c^\prime(0)=\bar v$. Then 
\begin{itemize}
\item $F_{H_X}$ is well defined and it is a norm on $T_{H_X}(G/H_X)$;
\item The function $F:T(G/H_X)\rightarrow \mathbb R$ defined by
\[
F(gH_X,\tilde v)=\lim_{t\rightarrow 0}\frac{d_X(c(t),gH_X)}{\vert t \vert},
\]
where $c:(-\varepsilon, \varepsilon) \rightarrow G/H_X$ is a curve such that $c(0)=gH_X$ and $c^\prime(0)=\tilde v$, does not depend on the choice of $c$. Moreover
\begin{equation}
\label{Fcontinuo}
F(gH_X,\tilde v):=F_{H_X}(d(\phi_{g^{-1}})_{gH_X}(\tilde v))
\end{equation}
holds for every $g \in gH_X$;
\item $F$ is a $G$-invariant Finsler metric;
\item $F_{H_X}$ is invariant by the isotropy representation $\iota: H_X \rightarrow \mathrm{Gl}(T_{H_X}(G/H_X))$, given by $\iota(h)=d(\phi_h)_{H_X}$, that is, $F(d(\phi_h)_{H_X}(v))=F(v)$ for every $(h,v)\in H_X \times T_{H_X}(G/H_X)$.
\end{itemize}
\end{theorem}

\

{\it Proof}

\

Item 1.

$F_{H_X}$ is well defined due to Corollary \ref{naodependedacurva}. 
Let us prove that $F_{H_X}$ is a norm.
The fact that $F_{H_X}\geq 0$ and $F_{H_X}(\bar 0)=0$ is straightforward from the definition of $F_{H_X}$. Observe that $F_{H_X}(\bar v)>0$ if $\bar v\neq 0$ is consequence of Theorem \ref{integralkillingdiferenciavel} and the fact that $\eta\vert_{[0,1]}>0$. 

Notice that $F_{H_X}(a\bar v)=\vert a\vert F_{H_X}(\bar v)$ for every $a\in \mathbb R$ and $\bar v\in T_{H_X}G/H_X$ is trivial when $a=0$. When $a\neq 0$, it is a direct consequence of Theorem \ref{integralkillingdiferenciavel}. 

Let us prove the triangle inequality.
Consider a decomposition $\mathfrak g=\mathfrak m \oplus \mathfrak h_X$.
Let $v,w\in \mathfrak m\subset \mathfrak g$ such that $v+\mathfrak h_X=\bar v$ and $w+\mathfrak h_X=\bar w$. Observe that
\[
d_X(H_X,\exp(tv+tw)H_X)
\]
\[
\leq d_X(H_X,\exp(tv)H_X)+d_X(\exp(tv)H_X,\exp(tu+tv)H_X)
\]
\[
= d_X(H_X,\exp(tv)H_X)+d_X(H_X,\exp(-tv) \exp(tu+tv) H_X)=C.
\]
But $\exp(-tv) \exp(tu+tv)=\exp(tu+O(t^2))$ (see for instance \cite{Helgas}). Then
\[
C\leq d_X(H_X,\exp(tv)H_X)+ d_X(H_X,\exp(tu)H_X) 
\]
\[
+d_X(\exp(tu+O(t^2))H_X,\exp(tu)H_X).
\]
\[
\leq d_X(H_X,\exp(tv)H_X)+ d_X(H_X,\exp(tu)H_X) 
\]
\[
+C_1 d_F(\exp(tu+O(t^2))H_X,\exp(tu)H_X)=C_2
\]
for some constant $C_1>0$, where $d_F$ is the distance function with respect to a Riemannian metric on $G/H_X$ (Theorem \ref{lipschitzlocalequivalente}).
Observe that $\exp:\mathfrak g\rightarrow G$ and $\pi:G\rightarrow G/H_X$ are smooth maps. 
If we endow $\mathfrak g$ with an arbitrary Euclidean metric $d_{\mathfrak g}$ then $\pi\circ \exp$ is locally Lipschitz. 
Then
\[
C_2\leq d_X(H_X,\exp(tv)H_X)+ d_X(H_X,\exp(tu)H_X) +C_3 d_{\mathfrak g}(tu+O(t^2),tu),
\]
and rewriting the inequality above we have that
\[
d_X(H_X,\exp(tv+tw)H_X) \leq d_X(H_X,\exp(tv)H_X)+ d_X(H_X,\exp(tu)H_X) + O(t^2).
\]
Now dividing the inequality above by $\vert t\vert$ and taking the limit with $t$ going to zero, we have the triangle inequality. 
Therefore $F_{H_X}$ is a norm on $T_{H_X}(G/H_X)$ and consequently it is continuous. 

\

\

Item 2.

Choose $g \in gH_X$ and consider $c:(-\varepsilon,\varepsilon) \rightarrow G/H_X$ such that $c(0)=gH_X$ and $c^\prime (0)=\tilde v$. Notice that
\[
F(gH_X,\tilde v)=\lim_{t\rightarrow 0}\frac{d_X(c(t),gH_X)}{\vert t \vert}=\lim_{t\rightarrow 0}\frac{d_X(g^{-1}c(t),H_X)}{\vert t \vert}
\]
\begin{equation}
\label{caracteriza F}
=F_{H_X}(d(\phi_{g^{-1}})_{gH_X}(\tilde v))
\end{equation}
due to the $G$-invariance of $d_X$. Observe that the equality 
\[
\lim_{t\rightarrow 0}\frac{d_X(c(t),gH_X)}{\vert t \vert}=F_{H_X}(d(\phi_{g^{-1}})_{gH_X}(\tilde v))
\] 
states that the left-hand side does not depend on the choice of $c$ such that $c(0)=gH_X$ and $c^\prime(0)=\tilde v$ and the right-hand side does not depend on the choice of $g\in gH_X$. Observe that (\ref{invariante continua}) applies here.

\

Item 3

If we suppose that we are in the settings \ref{mhx}, we have that $F$ is continuous in $T(B_{d_\mathfrak m}(H_X,r))$ due to (\ref{invariante continua}). 
The proof that $F$ is $G$-invariant is analogous to (\ref{caracteriza F}). 
Therefore $F$ is a $G$-invariant Finsler metric.

\

Item 4

It is a direct consequence of (\ref{caracteriza F}) with $g \in H_X$. $\blacksquare$

\

Suppose that we are in the hypotheses of Theorem \ref{intrinsecaefinsler}. We have that
\[
\hat d_X(x,y)=\inf_{\gamma \in \mathcal C_{x,y}}\ell_{d_X}(\gamma)
\]
but
\[
d_F(x,y)=\inf_{\gamma \in \mathcal S_{x,y}}\ell_{d_X}(\gamma).
\]
The next lemma implies that $\hat d_X=d_F$.


\begin{lemma}
\label{intrinsecaefinsler2}
Let $\varphi:G\times M\rightarrow M$ be a smooth action by isometries of a Lie group $G$ on a differentiable manifold $(M,d\in \mathcal L(M))$ and $X\subset M$ be a compact subset. Let $\gamma:[a,b]\rightarrow (G/H_X,d_X)$ be a rectifiable path. Then
\begin{itemize}
\item $\gamma$ admits a reparameterization  $\eta:[a^\prime,b^\prime]\rightarrow (G/H_X,d_X)$ such that $\eta^\prime$ exists a.e.. Moreover
\begin{equation}
\label{eqintrinsecafinsler}
\ell_{\hat d_X}(\gamma)=\ell_{d_X}(\gamma)=\ell_{d_X}(\eta)=\int_{a^\prime}^{b^\prime} v_{\eta}(t)dt=\int_{a^\prime}^{b^\prime} F(\eta(t),\eta^\prime(t))dt.
\end{equation}

\item For every $\varepsilon > 0$, then there exist a path $\gamma_\varepsilon:[a,b] \rightarrow (G/H_X,d_X)$ in $\mathcal S_{\gamma(a),\gamma(b)}$, such that $\ell_{d_X}(\gamma_\varepsilon) < \ell_{d_X}(\gamma)+\varepsilon$.

\end{itemize} 
\end{lemma}

Item 1.

Let $\gamma:[a,b] \rightarrow (G/H_X,d_X)$ be a rectifiable path. 
Then there exist a Lipschitz reparameterization $\eta:[a^\prime,b^\prime]\rightarrow (G/H_X,d_X)$ of $\gamma$ (see Proposition 2.5.9 of \cite{Burag}). 
Notice that for every $p\in G/H_X$, there exist a $\tau$-neighborhood $V$ of $p\in G/H_X$ and an Euclidean metric $d_E$ on $V$ such that $d_X\vert_V$ and $d_E$ are Lipschitz equivalent (Theorem \ref{lipschitzlocalequivalente}).
Let $(x_1,\ldots,x_n)$ be a coordinate system on $V$ such that $\{\partial/\partial x_1,\ldots, \partial /\partial x_n\}$ is an orthonormal basis of $(V, d_E)$.
If we restrict $\eta$ in such a way that its image is contained in $V$, then $\eta:[c^\prime,d^\prime]\rightarrow (V,d_E)$ is Lipschitz, as well as their coordinate functions $x_i \circ \eta$.
Then the derivative of every coordinate function exits a.e. what implies that the derivative of $\eta$ exists a.e.. 
Then the last equality of (\ref{eqintrinsecafinsler}) holds due to Theorem \ref{funcaoF}. The third equality follows from (\ref{integral velocidade escalar}) and the other equalities are immediate, what settle (\ref{eqintrinsecafinsler}).

\

Item 2.

Suppose that $\gamma$ is reparameterized as $\eta$ in Item 1 such that (\ref{eqintrinsecafinsler}) holds. Set $x=\gamma(a)$ and $y=\gamma(b)$. 
In order to prove the second item, suppose that 
\[
\inf_{\tilde \gamma\in \mathcal S_{x,y}}\ell_{d_X}(\tilde \gamma)
> \ell_{d_X}(\gamma).
\]
We claim for a contradiction. Consider
\[
C=\frac{\inf\limits_{\tilde \gamma\in \mathcal S_{x,y}}\ell_{d_X}(\tilde \gamma)}{\ell_{d_X}(\gamma)}
> 1
\]
and let $\varepsilon >0$ such that $(1+\varepsilon)^2< C$. Then for every $\tilde \delta >0$, there exist an interval $[a^\prime,b^\prime] \subset [a,b]$ with $\ell([a^\prime, b^\prime])<\tilde \delta$ such that
 
\begin{equation}
\label{suave maior}
\frac{\inf\limits_{\tilde \gamma\in \mathcal S_{\gamma(a^\prime),\gamma(b^\prime)}}\ell_{d_X}(\tilde \gamma)}{\ell_{d_X}(\gamma\vert_{[a^\prime,b^\prime]})} \geq C.
\end{equation}
In fact, if we subdivide $[a,b]$ in several subintervals, the restriction of $\gamma$ to one of these subintervals should satisfy (\ref{suave maior}).

For every $u\in [a,b]$, there exist a $\tau$-neighborhood $V_u$ of $\gamma(u)$ and a Finsler metric $F_u$ on $V_u$ such that 
\begin{itemize}
\item $(V_u,F_u)$ is isometric (as a Finsler manifold) to a convex open subset of $(\mathbb R^n,F_0)$; 
\item $F_u(\gamma(u),\cdot)=F(\gamma(u),\cdot)$;
\item $(1+\varepsilon)^{-1}F_u \leq F \leq (1+\varepsilon ) F_u$ (see Proposition \ref{vizinhancas lipschitz c perto de um}).
\end{itemize}
This induces an open cover of $[a,b]$ by intervals such that $\gamma$ restricted to these intervals have their image contained in some $V_u$.
Let $\delta>0$ be the Lebesgue number of this covering. 
Then we can find a closed interval $[\hat a, \hat b]$ of length less than $\delta>0$ such that (\ref{suave maior}) holds. But
\[
\ell_{d_X}(\gamma\vert_{[\hat a,\hat b]}) 
= \int_{\hat a}^{\hat b}F(\gamma(t),\gamma^\prime (t)) dt 
\geq (1+\varepsilon)^{-1} \int_{\hat a}^{\hat b} F_u(\gamma(t),\gamma^\prime(t))dt
\]
\[
\geq (1+\varepsilon)^{-1}\int_{\hat a}^{\hat b} F_u(\alpha(t),\alpha^\prime (t))dt \geq (1+\varepsilon)^{-2}\int_{\hat a}^{\hat b} F(\alpha(t),\alpha^\prime(t))dt
\]
\[
=(1+\varepsilon)^{-2}\ell_{d_X}(\alpha),
\]
where $\alpha$ is the straight line connecting $\gamma(\hat a)$ and $\gamma(\hat b)$ on $(V_u, F_u)\approx U\subset (\mathbb R^n,F_0)$, what contradicts (\ref{suave maior}). This settles the second item and the lemma. $\blacksquare$

\begin{theorem}
\label{intrinsecaefinsler}
Let $\varphi:G\times M\rightarrow M$ be a smooth action by isometries of a Lie group $G$ on a differentiable manifold $(M,d\in \mathcal L(M))$ and $X\subset M$ be a compact subset. If $\gamma:[a,b]\rightarrow G/H_X$ is a path which is continuously differentiable by parts (or more in general Lipschitz), then the speed $v_{\gamma}(t)$ exists whenever $\gamma^\prime(t)$ exists,
\begin{equation}
\label{eqintrinsecafinsler2}
\ell_{\hat d_X}(\gamma)=\ell_{d_X}(\gamma)=\int_a^b v_{\gamma}(t)dt=\int_a^b F(\gamma(t),\gamma^\prime(t))dt
\end{equation}
and $(G/H_X,\hat d_X)$ is Finsler (the integrals in (\ref{eqintrinsecafinsler2}) are Lebesgue integrals).
\end{theorem}

\

{\it Proof}

\

First of all, we prove (\ref{eqintrinsecafinsler2}) with $v_\gamma(t)$ existing a.e.. 
By Proposition \ref{suaveretificavel}, every path in $(G/H_X,d_X)$ which is continuously differentiable by parts is Lipschitz (and $d_X$-rectifiable). This settles the first equation of (\ref{eqintrinsecafinsler2}) and implies that $v_\gamma(t)$ exists a.e..
The second equation of (\ref{eqintrinsecafinsler2}) holds due to (\ref{integral velocidade escalar}). 
The last equation of (\ref{eqintrinsecafinsler2}) follows from the definition of $F$. 

In order to see that $(G/H_X,\hat d_X)$ is Finsler, notice that (\ref{eqintrinsecafinsler2}) implies that $d_X$ and $F$ induces the same length function on paths which are continuously differentiable by parts. Moreover
\[
\hat d_X(x,y)=\inf\limits_{\gamma \in \mathcal C_{x,y}}\ell_{d_X}(\gamma) =\inf\limits_{\gamma \in \mathcal S_{x,y}}\ell_{d_X}(\gamma)=d_F(x,y)
\]
where the second equality is due to the second item of Lemma \ref{intrinsecaefinsler2}. 
Then $(G/H_X,\hat d_X)$ is isometric to $(G/H_X,d_F)$ and $\hat d_X$ is Finsler.

Finally $v_\gamma(t)$ exists whenever $\gamma^\prime(t)$ exists due to Theorem \ref{velocidade Finsler}.$\blacksquare$

\begin{corollary}
\label{intrinsecafinsler2}
Let $\varphi:G\times M\rightarrow M$ be a transitive smooth action by isometries of a Lie group $G$ on a differentiable manifold $(M,d\in \mathcal L (M))$, where $d$ is an intrinsic metric. Then $(M,d)$ is Finsler.
\end{corollary}

{\it Proof}

\

Consider $p\in M$ and set $X=\{p\}$. Then $(M,d)$ is isometric to $(G/H_X,d_X)$. But $d_X=\hat d_X$ because $d$ is intrinsic. But $\hat d_X$ is Finsler due to Theorem \ref{intrinsecaefinsler}, what settles the corollary.$\blacksquare$

\section{Further examples}
\label{exemplos}

This section is devoted to present additional examples in order to better illustrate the theory developed so far.

\begin{example} (Counterexample to Corollary \ref{intrinsecafinsler2}) If a metric $d$ on a differentiable manifold $M$ is Lipschitz equivalent to a Finsler metric $F$ but it is not invariant by a transitive action of a group, we do not have necessarily that $d$ is Finsler: 
Let $f:\mathbb R\rightarrow \mathbb R$ be a strictly increasing Lipschitz function which is not differentiable at $0$. 
Moreover we can also suppose that the left and right derivatives does not exist at $0$ and that there exist constants $C_1,C_2>0$ such that $C_1\vert x-y\vert\leq f(x-y)\leq C_2 \vert x-y\vert$. 
For instance, consider the lines $L_1:=\{(x,y)\in \mathbb R^2;y=x\}$ and $L_2:=\{(x,y)\in \mathbb R^2;y=x/2\}$ in $\mathbb R^2$. 
The graph of $f$ will be placed between these lines. 
It is the concatenation of the line segments that connect the points $(1,1)\in L_1$, $(3/4,3/8)\in L_2$, $(1/4,1/4)\in L_1$, $(3/16,3/32)\in L_2$, $(1/16,1/16)\in L_1$ and so on. 
We proceed in a similar way for $t>1$ and set $f(t):=-f(-t)$ for $t<0$. 
Finally we define $f(0)=0$.
The function $f$ satisfies the conditions stated above.

Consider $\mathbb R$ with the metric $d(x,y)=\vert f(x)-f(y)\vert$. 
Of course $(\mathbb R,d)$ is intrinsic and Lipschitz equivalent to $\mathbb R$ with the canonical metric. If we suppose that the Finsler metric $F$ such that $d_F=d$ exists, we would have that
\[
F(v)=\lim_{t\rightarrow 0+}\frac{f(tv-0)}{t}
\]
which is not possible. Therefore $d$ is not Finsler.
\end{example}

\begin{example}
If $\varphi:G\times M\rightarrow M$ is a smooth transitive action by isometries of a Lie group on a differentiable manifold $(M,d\in \mathcal L (M))$, then we do not have necessarily that $d$ is intrinsic. 
In fact, we have the following example: $G=(\mathbb R,+)$, $M=\mathbb R$ is endowed with the metric
\[
d(x,y)=
\left\{
\begin{array}{ccc}
\vert x-y \vert & \text{ if } & \vert x-y\vert\leq 1;\\
1 & \text{ if } & \vert x-y\vert> 1,
\end{array}
\right.
\]
and $\varphi: G\times M\rightarrow M$ is the usual sum of real numbers. Observe that $d$ is not intrinsic (for instance, the length of $\gamma:[0,2]\rightarrow \mathbb R$ given by $\gamma(x)=x$ is equal to $2$).\end{example}

\begin{example}
\label{rnrn}
{\em (Translation in Euclidean space)} Let $G=(\mathbb R^n,+)$, $M=(\mathbb R^n,d)$, where $d$ is the Euclidean metric and consider $\varphi:\mathbb R^n\times \mathbb R^n\rightarrow \mathbb R^n$ be the action given by $\varphi(x,y)=x+y$.
Let $X\subset \mathbb R^n$ be a non empty compact set. 
The isotropy subgroup is the trivial subgroup of $\mathbb R^n$ what implies that $\mathbb R^n/H_X\cong\mathbb R^n$. 
The equalities  $d_X(xH_X,yH_X)=d_H(x+X,y+X)=d_H(x-y+X,X)$ hold because the action is by isometries. We claim that $d_X(xH_X,yH_X)=d(x,y)$.
 
The relationship $d_H(x+X,y+X)\leq d(x,y)$ is straightforward due to the definition of Hausdorff distance and because $\varphi$ is an action by isometries. 

In order to see that $d_H(x+X,y+X)\geq d(x,y)$, set $z=x-y$. 
By an orthogonal change of coordinate systems, we can suppose that $z=(0,\ldots, 0, z_n)$, $z_n\geq 0$.
Consider $ \pi_n:\mathbb R^n\rightarrow \mathbb R$ the projection in the $n$-th coordinate, set $m=\min_{x\in X}\pi_n(x)$ and consider $\bar x\in X$ such that $\pi_n(\bar x)=m$. 
Notice that $d_H(X,z+X)\geq d(\bar x,z+X)\geq d(\bar x, \bar x+z)$ due to the definition of Hausdorff distance and the choice of $\bar x$. This settles $d_H(x+X,y+X)\geq d(x,y)$.
Therefore $d_H(x+X,y+X)=d(x,y)$ and $id:(\mathbb R^n,d_X) \rightarrow (\mathbb R^n,d)$ is an isometry. Observe that $\hat d_X=d_X$ because $d_X=d$ is an intrinsic metric. 
\end{example}

Before the next example, we prove the following lemma.

\begin{lemma}
\label{intrinsecaigual}
Let $\varphi:G\times M\rightarrow M$ be a smooth transitive action of a Lie group on a differentiable manifold $M$, $p\in M$ and suppose that $d$ and $\rho$ are two $G$-invariant metrics on $M$ such that $d(p,gp)=\rho(p,gp)$ for every $g$ in a neighborhood $V\subset G$ of $e$. 
Then $\hat d=\hat \rho$.
\end{lemma}

{\it Proof}

\

Let $U$ be a neighborhood of $e\in G$ such that $U^2\subset V$ and $U=U^{-1}$. 
Let $\eta:[a,b]\rightarrow M$ be a path on $M$. 
For every $t\in [a,b]$, choose $g_t\in G$ such that $g_t\eta(t)=p$. 
Consider the open cover $\mathcal O=\{\eta^{-1}(g_t^{-1}Up);t\in [a,b]\}$ of $[a,b]$ and denote by $\delta>0$ a Lebesgue number of $\mathcal O$.
Fix a partition $\mathcal P=\{a=t_0<t_1<\ldots < t_{N_\mathcal P}=b\}$ such that $\vert \mathcal P\vert<\delta$.
For every $i=1,\ldots,N_{\mathcal P}$, consider $\tilde t_i\in [a,b]$ such that $g_{\tilde t_i}\eta([t_{i-1},t_i])\subset Up$.
Then there exist $h_{t_i},u_{t_i}\in U$ such that $g_{\tilde t_i}\eta(t_{i-1})=h_{t_i} p$ and $g_{\tilde t_i}\eta(t_i)=u_{t_i} p$. Observe that $h_{t_i}^{-1}u_{t_i} \in V$. Then

\[
\sum_{i=1}^{N_{\mathcal P}}d(\eta(t_{i-1}),\eta(t_i))=\sum_{i=1}^{N_{\mathcal P}}d(g^{-1}_{\tilde t_i}h_{t_i}p,g^{-1}_{\tilde t_i}u_{t_i}p)
=\sum_{i=1}^{N_{\mathcal P}}d(p,h_{t_i}^{-1}u_{t_i}p)
\]
\[
=\sum_{i=1}^{N_{\mathcal P}}\rho(p,h_{t_i}^{-1}u_{t_i}p)=\sum_{i=1}^{N_{\mathcal P}}\rho(\eta(t_{i-1}),\eta(t_i)).
\]
Taking the supremum with respect to $\mathcal P$, we have that $\ell_d(\eta)=\ell_\rho(\eta)$ for every path $\eta$. Thus $\hat d=\hat \rho$.$\blacksquare$

\begin{example}
\label{bolageodesica}
Let $\varphi:G\times M\rightarrow M$ a transitive action by isometries of a Lie group $G$ on a Riemannian manifold $M$. Denote the Riemannian distance by $d$.
Let $X$ be the closure of geodesic ball in $M$. 
We remember that $B(p,r)$ is a geodesic ball if its closure is contained in a normal neighborhood of $p$ (see \cite{Carmo3}).
Notice that $H_X=H_p$ because an isometry sends a geodesic ball of radius $r$ in a geodesic ball of the same radius and $hX=X$ if and only if $hp=p$. Thus $G/H_X$ is diffeomorphic to $M$.

Due to the definition of geodesic ball, there exist an $\varepsilon>0$ such that $B(p,r+\varepsilon)$ is also a geodesic ball. 
Then there exist a neighborhood $V$ of $e\in G$ such that $V=V^{-1}$ and $B(gp,r)=gB(p,r)\subset B(p,r+\varepsilon)$ for every $g\in V$. For $g\in V-H_p$, consider the unique geodesic $\gamma_g$ such that $\gamma_g(0)=p$ and $\gamma_g(1)=gp$. 
Let $\bar t$ be the first positive value such that $\gamma_g(\bar t)\in \partial B(p,r)$.
Then $d(p,\gamma_g(\bar t))=d(gp,\gamma_g(\bar t+1))=r$ and $d(\gamma_g(\bar t),\gamma_g(\bar t+1))=d(p,gp)$, what implies that $\gamma_g(\bar t+1)\in \partial B(p,r+d(p,gp))$.
Observe that $p$, $gp$, $\gamma_g(\bar t)$ and $\gamma_g(\bar t+1)$ are in the same minimizing geodesic.

In order to prove that $gH_X\mapsto gp$ is an isometry from $(G/H_X,\hat d_X)$ to $(M,d)$, it is enough to prove that $d_H(X,gX)=d(p,gp)$ for every $g\in V$ due to Lemma \ref{intrinsecaigual}.
If $g\in H_p=H_X$, there is nothing to prove.
Otherwise notice that $d_H(X,gX)\geq d(X,\gamma_g(\bar t+1))$ holds due to the definition of $d_H$ and the fact that $\gamma_g(\bar t+1)\in gX$.
Moreover the inequality
\[
d(x,\gamma_g(\bar t+1))\geq d(p,\gamma_g(\bar t+1))-d(p,x)
\]
holds for every $x\in X$. If we consider the infimum of the right hand side with respect to $x\in X$, we get
\[
d(x,\gamma_g(\bar t+1))\geq d(p,\gamma_g(\bar t+1))-d(p,\gamma_g(\bar t))
\]
\[
=d(\gamma_g(\bar t),\gamma_g(\bar t+1))= d(p,gp) \text{ for every }x\in X.
\]
Therefore $d(X,\gamma_g(\bar t+1))\geq d(p,gp)$ and $d_H(X,gX)\geq d(p,gp)$ for every $g\in V$.

In order to see that $d(p,gp)\geq d_H(X,gX)$ we use (\ref{hausdorffalternativo}). 
It is straightforward that $gX=\bar B(gp,r)\subset \bar B(p,r+d(p,gp))=\bar B(X,d(p,gp))$. 
In order to see that $X\subset \bar B(gX,d(p,gp))$, observe that $\varphi(g,\cdot):M\rightarrow M$ is an isometry that sends $B(p,r+d(p,q))$ to $B(gp,r+d(p,q))$ and both are geodesic balls. 
Then $\bar B(gp,r+d(p,q))=\bar B(gX,d(p,gp))$ and $X\subset \bar B(gX,d(p,gp))$ holds.
Consequently $d(p,gp)\geq d_H(X,gX)$.

\end{example}

\

Example \ref{rnrn} shows a situation where $(G/H_X,d_X)$ and $(G/H_X,\hat d_X)$ are isometric to $(M,d)$. 
In Example \ref{bolageodesica} we have that $(G/H_X,\hat d_X)$ is isometric to $(M,d)$ but $(G/H_X,d_X)$ is not necessarily isometric to $(M,d)$. 
We also have examples where $(G/H_X,\hat d_X)$ is not isometric to $(M,d)$, even if $G/H_X$ is diffeomorphic to $M$. 
Example \ref{toroplano} shows a compact case and Example \ref{hiperbolico} shows a noncompact case.

\begin{theorem}
\label{subvariedade}
Let $(M,\left<\cdot,\cdot\right>)$ be a Riemannian manifold and $\varphi:G\times M\rightarrow M$ be a smooth action by isometries of a Lie group $G$ on $M$. 
Let $X\subset M$ be a compact submanifold without boundary. 
Then the Finsler metric $F$ induced by the metric $\hat d_X$ on $G/H_X$ is given by
\begin{equation}
\label{Fmergulhado}
F(v+\mathfrak h_X)=\max_{x\in X}\left\Vert \left( K_v(x) \right)^N\right\Vert_M
\end{equation}
for every $v\in \mathfrak g$.
\end{theorem}

Before proving Theorem \ref{subvariedade}, we prove two lemmas.

\begin{lemma}
\label{gauss}
Let $(V_1,\left<\cdot,\cdot\right>_1)$ and $(V_2,\left<\cdot,\cdot\right>_2)$ be finite dimensional real vector spaces.
Consider orthogonal decompositions $V_1=T_1\oplus N_1$ and $V_2=T_2\oplus N_2$.  Let $\eta:V_1\rightarrow V_2$ be a linear map such that $\eta(T_1)\subset T_2$, $\eta(N_1)\subset N_2$ and such that $\eta\vert_{N_1}$ is an isometry over its image. 
Then $\left<\eta(v),\eta(w)\right>_2=\left<v,w\right>_1$ for every $(v,w)\in V_1\times N_1$.
\end{lemma}

{\it Proof}

\

Let $v\in V_1$ and consider the decomposition $v=v_{T_1}+v_{N_1}$, where $v_{T_1}$ and $v_{N_1}$ are the components of $v$ in $T_1$ and $N_1$ respectively.
Consider $w\in N_1$.
Then
\[
\left<\eta(v),\eta(w)\right>_2=\left<\eta(v_{T_1})+\eta(v_{N_1}),\eta(w)\right>_2=\left<\eta(v_{N_1}),\eta(w)\right>_2
\]
\[
=\left<v_{N_1},w\right>_1=\left<v_{T_1}+v_{N_1},w\right>_1
=\left<v,w\right>_1.\blacksquare
\]

\

The next lemma is interesting by itself because it does not impose any restriction on the compact subset $X$.

\begin{lemma}
\label{supcontinuo} 
Let $\varphi:G\times M\rightarrow M$ be a smooth action by isometries of a Lie group $G$ on a differentiable manifold $M$ endowed with a metric $d\in \mathcal L(M)$. 
Let $X$ be a compact subset of $M$. 
Suppose that the function $\bar f:X\times (\mathbb R-\{0\})\rightarrow \mathbb R$ given by $\bar f(x,t)=d_M(X,\exp(tv)x)/\vert t \vert$ can be continuously extended to a continuous function $f:X\times \mathbb R\rightarrow \mathbb R$. Then the Finsler metric correspondent to $\hat d_X$ is given by
\begin{equation}
\label{equacaoestendivel}
F(v+\mathfrak h_X)=\lim_{t\rightarrow 0}\frac{d_X(H_X,\exp(tv)H_X)}{\vert t\vert}=\sup_{x\in X} f(x,0).
\end{equation}
\end{lemma}

{\it Proof}

\

The first equality is due to Theorems \ref{funcaoF} and \ref{intrinsecaefinsler}.

For the sake of simplicity, denote $Y:=X\times \mathbb R$ and define $d_Y:Y\times Y\rightarrow \mathbb R$ by $d_Y((x_1,t_1),(x_2,t_2))=d_M(x_1,x_2)+\vert t_1-t_2\vert$. 
If $\varepsilon>0$ is given, then for every $(x,0)\in X\times \mathbb R$, there exist a $\delta_x>0$  such that $\vert f(y,t)- f(x,0)\vert<\varepsilon/2$ whenever $d_Y((y,t),(x,0))<\delta_x$. For every $x\in X$, define
\[
K(x,\delta_x):= \{(y,t)\in X\times \mathbb R; d_M(y,x)<\delta_x/2,\vert t \vert<\delta_x/2\},
\]
consider a finite open cover $\{K(x_i,\delta_{x_i});i=1,\ldots,k\}$ of $X\times \{0\}$ and set $\delta:=\min_{i=1,\ldots,k}\delta_{x_i}/2$. If $y\in X$, then there exist a $j\in \{1,\ldots, k\}$ such that $(y,0)\in K(x_j,\delta_{x_j})$. If $\vert t \vert<\delta$, then
\[
\vert f(y,t)-f(y,0)\vert\leq \vert f(y,t)-f(x_j,0)\vert+ \vert f(x_j,0)-f(y,0) \vert< \varepsilon. 
\]
In other words, $\vert f(y,t)-f(y,0)\vert <\varepsilon $ if $\vert t\vert<\delta$ for every $y\in X$.

Now fix $t\in (-\delta,\delta)$ and consider $\bar x,\hat x\in X$ such that $f(\hat x,t)=\sup_{x\in X}f(x,t)$ and $f(\bar x,0)=\sup_{x\in X}f(x,0)$. Then
\[
f(\bar x,0)- \varepsilon<f(\bar x,t)\leq f(\hat x,t)<f(\hat x,0)+\varepsilon\leq f(\bar x,0)+\varepsilon
\]
what implies that
\[
\vert \sup_{x\in X} f(x,t)-\sup_{x\in X}f(x,0) \vert < \varepsilon
\]
whenever $\vert t\vert<\delta$. Likewise
\[
\vert \sup_{x\in X} f(x,-t)-\sup_{x\in X}f(x,0) \vert < \varepsilon
\]
whenever $\vert t\vert<\delta$ and
\[
\vert \max\{\sup_{x\in X}f(x,t),\sup_{x\in X} f(x,-t)\}-\sup_{x\in X}f(x,0) \vert < \varepsilon
\]
whenever $\vert t\vert<\delta$. Then
\[
\lim_{t\rightarrow 0}\max\left\{\sup_{x\in X}\frac{d_M(X,\exp(tv)x)}{\vert t\vert},\sup_{x\in X}\frac{d_M(X,\exp(-tv)x)}{\vert t\vert}\right\}=\sup_{x\in X}f(x,0),
\]
what is the second equation of (\ref{equacaoestendivel}).$\blacksquare$

\

{\it Proof of Theorem \ref{subvariedade}}

\

Let $X$ be a compact submanifold of $M$ and denote the normal space of $X$ at $x$ by $T^\perp_x X$. Then we have the decomposition $T_xM=T_xX\oplus T^\perp_x X$. Consider the tubular neighborhood
\[
\{y\in M;y=\widetilde \exp_x (v), x\in X, v\in T_x^\perp X, \Vert v \Vert_M<\varepsilon \}
\]
of $X$ in $M$, where $\widetilde \exp$ is the Riemannian exponential and $\varepsilon>0$ is a sufficiently small positive value (see \cite{Spivak1}).

Using Theorems \ref{funcaoF} and \ref{intrinsecaefinsler} and the definition of Hausdorff distance, we have that
\[
F(v+\mathfrak h_X)=\lim_{t\rightarrow 0}\max\left\{\sup_{x\in X}\frac{d_M(x,\exp (tv)X)}{\vert t\vert} , \sup_{x\in X}\frac{d_M(X,\exp (-tv)x)}{\vert t\vert} \right\}.
\]

If we prove that the function $f:X\times \mathbb R\rightarrow \mathbb R$ given by
\begin{equation}
\label{velocidademedia}
f(x,t)=
\left\{
\begin{array}{ccc}
\frac{d_M(X,\exp(tv)x)}{\vert t\vert} & \text{ if } & t\neq 0; \\
\left\Vert \left( \left.\frac{d}{dt}\right\vert_{t=0} \exp (tv)x \right)^N\right\Vert_M & \text{ if } & t=0,
\end{array}
\right.
\end{equation}
is continuous, then the theorem is proved due to Lemma \ref{supcontinuo}. 

Let $m=\dim X$ and $m+n=\dim M$. Fix $x\in X$ and consider a coordinate system $\psi=(\psi_1,\ldots,\psi_m)$ on a neighborhood $W\subset X$ of $x$. 
In a eventually smaller neighborhood $W^\prime\subset W$ of $x$, define an orthonormal frame $(E_{m+1},\ldots, E_{m+n})$ of the normal bundle of $X\subset M$. Consider the subset 
\[
V:=\{y\in M;y=\widetilde \exp_x (v); x\in W^\prime, v\in T^\perp_xX, \Vert v\Vert_M<\varepsilon\}
\]
of the tubular neighborhood and define a coordinate system $(V,\Psi)$ given by 
\[
\Psi(\widetilde\exp_x(v))=(\psi_1(x),\ldots,\psi_m(x),a_{m+1}(x),\ldots, a_{m+n}(x)),
\]
where $v=\sum\limits_{i=m+1}^{m+n} a_i(v) E_i$.
Then there exist an $\varepsilon^\prime>0$ and a neighborhood $W^{\prime\prime}\subset W^\prime$ of $x$ such that the function $f\vert_{W^{\prime\prime}\times (-\varepsilon^\prime, \varepsilon^\prime)}$ written in terms of coordinate system $\Psi$ is given by
\begin{equation}
\label{velocidademediacoord}
g(x,t)=
\left\{
\begin{array}{ccc}
\sqrt{\sum_{j=1}^n \left( \frac{a_{m+j}(\exp(tv)x)}{t} \right)^2} & \text{ if } & t\neq 0; \\
\sqrt{\sum_{j=1}^n \left< \left. \frac{d}{dt} \right\vert_{t=0} \exp(tv)x, E_{m+j} \right>_M^2} & \text{ if } & t=0
\end{array}
\right.
\end{equation}
due to the definition of $\Psi$.

Consider the smooth function $\sigma_j:W^{\prime\prime} \times (-\varepsilon^\prime,\varepsilon^\prime)\rightarrow \mathbb R$ given by $\sigma_j(x,t)= a_{m+j}(\exp(tv)x)$. 
Notice that $\sigma_j(x,0)=0$ for every $x\in W^{\prime\prime}$, what implies that there exist a continuous function $\xi_j:W^{\prime\prime} \times (-\varepsilon^\prime,\varepsilon^\prime)\rightarrow \mathbb R$ such that $\sigma_j(x,t)=t\xi_j(x,t)$ and $\xi_j(x,0)=\left.\frac{\partial}{\partial t}\right\vert_{t=0}\sigma_j(x,t)$ (see \cite{Carmo3} or \cite{Kobay1}).
We claim that $g(x,t)=\sqrt{\sum_{i=1}^n \xi_j^2(x,t)}$ what settles the theorem. For $t\neq 0$, the equality is immediate. For $t=0$ we have that
\[
\xi_j(x,0)=\left.\frac{\partial}{\partial t}\right\vert_{t=0}a_{m+j}(\exp(tv)x)=\left.\frac{\partial}{\partial t}\right\vert_{t=0}\left<\Psi(\exp(tv)x),d\Psi_{(x,0)}(E_{m+j})\right>_{\mathbb R^n}
\]
\[
=\left<d\Psi_{(x,0)}\left(\left.\frac{\partial}{\partial t}\right\vert_{t=0}\exp(tv)x\right),d\Psi_{(x,0)}(E_{m+j})\right>_{\mathbb R^n}
\]
\[
=\left<\left.\frac{\partial}{\partial t}\right\vert_{t=0}\exp(tv)x,E_{m+j}\right>_{M}
\]
where the last equality is because $d\Psi_{(x,0)}:T_xX\oplus T^\perp_xX \rightarrow \mathbb R^m \oplus \mathbb R^n$ satisfies the conditions of Lemma \ref{gauss}, with $T_x^\perp X=N_1$. Then $g(z,t)=\sqrt{\sum_{i=1}^n \xi_j^2(z,t)}$ and $g$ is continuous.$\blacksquare$

\

A particular case of Theorem \ref{subvariedade} is when $X$ is a finite subset of $M$. 

\begin{corollary}
\label{finito}
Let $(M,\left<\cdot,\cdot\right>)$ be a Riemannian manifold and $\varphi:G\times M\rightarrow M$ be a smooth action by isometries of a Lie group $G$ on $M$. 
Let $X=\{x_1,\ldots,x_k\}\subset M$ a finite subset. 
Then the Finsler metric $F$ correspondent to the metric $\hat d_X$ on $G/H_X$ is given by
\begin{equation}
\label{Fmergulhado}
F(v+\mathfrak h_X)=\max_{i=1,\ldots,k}\left\Vert \left( K_v(x_i) \right)\right\Vert_M
\end{equation}
for every $v\in \mathfrak g$.
\end{corollary}

\begin{proposition}
\label{biinvariante}
Let $\varphi:G\times G\rightarrow G$ be the usual group product of a connected Lie group $G$, with $M=G$ endowed with a bi-invariant Riemannian metric. Let $X\subset M$ be a compact submanifold without boundary, $v\in \mathfrak g$ and $F$ be the Finsler metric correspondent to $(G/H_X,\hat d_X)$. Then the following statements are equivalent:
\begin{enumerate}
\item There exist $x\in X$ such that $K_v(x)\in T^\perp_xX$;
\item $F_{H_X}(v+\mathfrak h_X)=\Vert v\Vert_M$.
\end{enumerate}
\end{proposition}

{\it Proof}

\

If $\left<\cdot,\cdot \right>_M$ is bi-invariant, then
\[
\Vert K_v(g)\Vert_M=\Vert d(R_g)_e(v)\Vert_M=\Vert v \Vert_M
\]
for every $g\in G$, where the first equality is due to (\ref{killingemgrupos}).
In particular, its norm is constant along $M$. Then
\[
F_{H_X}(v+\mathfrak h_X)=\max_{x\in X}\left\Vert \left( K_v(x) \right)^N\right\Vert_M\leq \Vert v\Vert_M
\]
and the equality holds if and only if there exist a $x\in X$ such that $K_v(x)\in T^\perp_xX$.$\blacksquare$

\begin{example}
\label{toroplano} Let $G$ be the additive group $(\mathbb R^2,+)$ and $M=\mathbb R^2/\mathbb Z^2$ be the flat torus. 
We represent a point of $M$ by $(\bar x, \bar y)$, where $\bar x$ is the equivalence class of $x\in \mathbb R$ in $\mathbb R/\mathbb Z$. 
Consider the natural action $\varphi:\mathbb R^2\times M\rightarrow M$ defined by $\varphi(((g_1,g_2),(\bar x_1,\bar x_2))=(\overline{g_1+x_1},\overline{g_2+x_2})$ and set $X=M-Q$, where $Q$ is the equivalence class of the square $(1/4,3/4)\times (1/4,3/4)\subset \mathbb R^2$ in $M$. 
In order to make calculations, we represent $M$ by the square $[0,1]\times [0,1]\subset \mathbb R^2$ with its opposite sides identified.
Observe that $H_X=\mathbb Z\times \mathbb Z$ and that $G/H_X$ is diffeomorphic to $M$.
Moreover if $(g_1,g_2)\in (-1/10,1/10)\times (-1/10,1/10)\subset G$, then it is not difficult to see that $d_X((g_1,g_2)H_X,H_X)=\max\{\vert g_1\vert,\vert g_2\vert\}$. 
Now we can use Lemma \ref{intrinsecaigual} to conclude that $(G/H_X,\hat d_X)$ is isometric to $(M,F)$, where $F$ is the homogeneous Finsler metric such that $F_{H_X}((v_1,v_2))=\max\{\vert v_1\vert, \vert v_2\vert \}$. 
Observe that $G/H_X$ is diffeomorphic to $M$ but $(G/H_X,\hat d_X)$ and $(M,d)$ are not isometric.
\end{example}

\begin{example}
\label{hiperbolico}
Let $M=\mathbb H^2=\{(x_1,x_2)\in \mathbb R^2;x_2>0\}$ be Poincar\'e half-plane model of the hyperbolic plane. 
Consider the Lie group $G=\{g_1,g_2)\in \mathbb R^2;g_2>0\}$ with the product $(g_1,g_2).(h_1,h_2)=(g_2h_1+g_1,g_2h_2)$ and inverse $(g_1,g_2)^{-1}=(-g_1/g_2,1/g_2)$.
Observe that $G$ is the semidirect product $\mathbb R \rtimes_\zeta \mathbb R^+$ of the additive group of real numbers and the multiplicative group of strictly positive real numbers, where the homomorphism $\zeta:(\mathbb R^+,0)\rightarrow \text{Aut}(\mathbb R)$ is such that $\zeta(t)$ is the multiplication by $t$.

Consider the action of $G$ on $M$ given by $\varphi:G\times M\rightarrow M$, $\varphi((g_1,g_2),(x_1,x_2))=g_2(x_1,x_2)+(g_1,0)$. 
Observe that $\varphi_{(g_1,g_2)}$ is a homothety $(x_1,x_2)\mapsto g_2(x_1,x_2)$ followed by a translation $(x_1,x_2)\mapsto (g_1+x_1,x_2)$, what implies that 
$\varphi$ is a smooth transitive action of $G$ on $M$ by isometries. 
Consider $X=\{(-a,b),(a,b)\}\subset M$, with $a,b>0$. 
Notice that $H_X=\{e\}$ and $G/H_X$ is diffeomorphic to $M$. 

Now we determine the Finsler function $F$ correspondent to $\hat d_X$. 
It is enough to determine $F_{H_X}(\cos \theta,\sin \theta)$ for $\theta\in [0,2\pi)$.
The point $e=(0,1)\cong H_X$ is the identity element of $G$. Then we can consider curves $c_\theta(t)=(t\cos\theta,t\sin\theta+1)$ and
\[
F_{H_X}((\cos\theta,\sin\theta))=\max\left\{ \left\Vert \left.\frac{d}{dt} \right\vert_{t=0} c_\theta(t)(-a,b) \right\Vert, \left\Vert \left.\frac{d}{dt} \right\vert_{t=0} c_\theta(t)(a,b) \right\Vert \right\}
\]
due to Theorems \ref{funcaoF}, \ref{intrinsecaefinsler} and Corollary \ref{finito}.
Direct calculations show that
\[
\left.\frac{d}{dt}\right\vert_{t=0}c_\theta(t)(x,y)=(x\sin\theta+\cos\theta,y\sin\theta)
\]
and
\[
\left\Vert \left.\frac{d}{dt}\right\vert_{t=0}c_\theta(t)(x,y)\right\Vert_{\mathbb H^2}= \frac{1}{y}\sqrt{x\sin2\theta +\cos^2\theta+(x^2+y^2)\sin^2\theta}.
\]
If we denote $v_\theta=(\cos\theta,\sin\theta)$, then
\[
F_{H_X}(v_\theta)=
\left\{
\begin{array}{ccc}
\frac{1}{b}\sqrt{a\sin2\theta +\cos^2\theta+(a^2+b^2)\sin^2\theta} & \text{ if } & \sin 2\theta >0 \\
\frac{1}{b}\sqrt{-a\sin2\theta +\cos^2\theta+(a^2+b^2)\sin^2\theta} &  & \text{otherwise}.
\end{array}
\right.
\]
Straightforward calculations show that $F$ is not smooth.
\end{example}

\end{document}